\documentclass[12pt,a4paper]{article}
\usepackage[truedimen,margin=30mm]{geometry} 

\usepackage{mathrsfs}
\usepackage{amssymb}
\usepackage{amsmath}
\usepackage{ascmac}
\usepackage{amsthm}
\usepackage[dvips]{graphicx}
\usepackage{natbib}
\usepackage{setspace}
\usepackage{times}
\usepackage{comment}
\usepackage{color}
\usepackage{multicol, multirow}
\usepackage{lscape}
\usepackage{natbib}
\usepackage{booktabs}
\usepackage{siunitx}
\usepackage[colorlinks=true,
            citecolor=blue]{hyperref}

\usepackage{titlesec}
\titleformat*{\section}{\large\bfseries}
\titleformat*{\subsection}{\it}

\newtheorem{thm}{Theorem}

\newtheorem{cor}{Corollary}
\newtheorem{prp}{Proposition}

\newtheorem{algo}{Algorithm}

\def\al{{\alpha}}

\def\ep{{\varepsilon}}

\def\bal{{\text{\boldmath $\alpha$}}}
\def\bbe{{\text{\boldmath $\beta$}}}
\def\bga{{\text{\boldmath $\gamma$}}}
\def\bta{{\text{\boldmath $\eta$}}}

\def\balh{{\text{\boldmath $\widehat\alpha$}}}
\def\bbeh{{\text{\boldmath $\widehat\beta$}}}
\def\bgah{{\text{\boldmath $\widehat\gamma$}}}

\def\Var{{\rm Var}}

\def\D{{\text{\boldmath $D$}}}
\def\W{{\text{\boldmath $W$}}}
\def\X{{\text{\boldmath $X$}}}
\def\Z{{\text{\boldmath $Z$}}}

\def\x{{\text{\boldmath $x$}}}
\def\z{{\text{\boldmath $z$}}}

\def\one{{\bf\text{\boldmath $1$}}}

\title{{\bf Causal Small Area Estimation with \\
Survey-only Covariates}}

\date{}
\author{}

\begin{document}

\maketitle
\doublespacing

\vspace{-1.5cm}
\begin{center}
{\large 
Tsubasa Ito$^1$ and Shonosuke Sugasawa$^2$
}

\medskip
\noindent
$^1$Faculty of Economics and Business, Hokkaido University\\
$^2$Faculty of Economics, Keio University
\end{center}

\medskip
\begin{center}
{\bf \large Abstract}
\end{center}

\vspace{-0cm}

\noindent
Area-specific causal inference is important in many policy and survey applications, where the goal is to evaluate treatment effects for small geographic or demographic domains. Existing causal small area estimation methods, however, typically rely on a strong data requirement that treatment status is observed for all units in the population. This assumption is often unrealistic in practical survey settings, where both treatment and outcome variables are observed only for sampled units, while auxiliary covariates are available for the full population. To address this limitation, we develop a new identification strategy for area-specific treatment effects under this more realistic data structure by combining survey-only covariates with population-level auxiliary information. Based on this result, we propose a doubly robust estimator that remains consistent when either the outcome regression model or the treatment and area assignment models are correctly specified. We further derive the semiparametric efficiency bound for the target parameter and show that the proposed estimator attains this bound under regularity conditions. Simulation studies demonstrate favorable finite-sample performance, particularly in settings with small sample sizes within areas, and an empirical application illustrates the practical relevance of the proposed framework.

\vspace{4mm}
\noindent
{\it Key words and phrases:}
Causal inference; Small area estimation; Survey data; Identification; Doubly robust estimation; Semiparametric efficiency

\section{Introduction}\label{sec:int}

Estimating population quantities for small geographic or demographic areas is a central task in official statistics, with important applications in policy evaluation and resource allocation. 
Because sample sizes within each area are often limited, direct estimators can be highly variable. To address this issue, the small area estimation (SAE) literature has developed methods that borrow strength across areas by incorporating auxiliary information available at the population level, such as census data \citep{RI15, F79, ELL03, BHF88, GR94, sugasawa2020small}. 
In many applications, however, the primary objective is not merely descriptive but causal, requiring estimation of area-specific treatment effects. 
This has led to recent efforts to integrate SAE methodology with causal inference.

Existing causal SAE methods, such as \citet{RSP23} and \citet{R25}, combine mixed-effects models with imputation of unobserved potential outcomes. 
These methods extend classical SAE frameworks by incorporating treatment variables and estimating area-specific average treatment effects through model-based predictions. 
A key assumption underlying these approaches is that the treatment variable is observed for all units in the population. 
Under this assumption, missing potential outcomes can be predicted using outcome regression models, and area-level causal effects are constructed by aggregating predicted values. 
However, this assumption is often unrealistic in practice. 
In many survey settings, both the outcome and the treatment are observed only for sampled individuals, while auxiliary covariates are available for the full population through administrative or census data \citep{S03, PS99}. 
This mismatch between the assumptions underlying existing methods and realistic data structures poses a fundamental challenge for causal small area estimation. 
When treatment status itself is only partially observed, the imputation-based strategies proposed in the current SAE literature \citep{RSP23, R25} are not directly applicable. 
Furthermore, these model-based approaches rely heavily on correct specification of outcome regression models. In complex and heterogeneous populations, model misspecification is a serious concern, and such methods do not possess the doubly robust (DR) property that is central in modern causal inference \citep{RRZ94, BR05, CCD18, H98, T06}. 
As a result, their theoretical guarantees are limited and their performance may deteriorate when modeling assumptions are violated.

This discrepancy between data sources highlights the need for methods that can effectively integrate survey and population information. 
To address this challenge, in this paper, we propose a new framework for causal small area estimation that explicitly accounts for the fact that both treatment and outcome are observed only for sampled units. 
Our approach integrates design-based inference \citep{HT52, H71, S03, B83} with doubly robust estimation, and crucially exploits rich covariates observed only for sampled units. 
These survey-only variables play a central role in enabling identification of small area treatment effects, as they provide information needed to adjust for systematic differences across areas in both treatment assignment and outcome-related characteristics that are not captured by population-level covariates alone. 
By leveraging this additional information, the proposed framework makes it possible to identify causal effects that are otherwise not recoverable under conventional small area estimation approaches. Our approach is also closely related to the identification strategy developed in \citet{KY21}, which derives an identification formula for small area means based on inverse probability weighting. 
While their framework focuses on point estimation of population averages, it does not explicitly consider causal parameters or treatment effects. We extend this line of work by embedding the identification strategy into the potential outcomes framework \citep{R74, RR83} and developing estimators for small area treatment effects under survey sampling.

A key feature of the proposed estimator is its robustness to model misspecification, namely, it remains consistent even when not all nuisance functions are correctly specified. 
In particular, this doubly robust structure arises naturally from the above identification framework, and provides both robustness against model misspecification and improved efficiency. 
We further derive the semiparametric efficiency bound for the target parameter and show that the proposed estimator attains this bound under regularity conditions. 
These results build on the semiparametric theory of causal inference \citep{RRZ94, N94, V00} and extend it to the small area estimation setting under complex sampling. From a practical standpoint, the proposed framework enables stable estimation of area-specific treatment effects even when sample sizes are small. 
By borrowing information across areas and adjusting for selection bias through propensity weighting, the method improves efficiency relative to direct estimators while mitigating the instability of inverse probability weighting estimators \citep{HT52, HI03}. 
At the same time, it avoids the strong modeling assumptions required by mixed-effects imputation approaches \citep{RSP23, R25}. 
We investigate the finite-sample performance of the proposed method through extensive simulation studies. 
The results demonstrate substantial efficiency gains over direct estimators and improved finite-sample stability relative to standard inverse probability weighting approaches. 
We also examine the accuracy of variance estimators and highlight challenges in variance estimation in finite samples. An application to real data further illustrates the practical relevance of the proposed approach. 
Overall, this paper contributes to the literature by introducing a new design-based and doubly robust framework for causal small area estimation under partial observability of treatment. 
By bridging the gap between causal inference and small area estimation, our approach provides a flexible and theoretically grounded tool for policy-relevant analysis in modern data environments.

The remainder of this paper is organized as follows.
Section \ref{sec:set} describes the setup, target parameters and identification of area-specific average treatment effects.
Section \ref{sec:est} discusses the causal small area estimators and their computation algorithm with inference strategy.
Section \ref{sec:sim} presents simulation results.
Section \ref{sec:rda} illustrates the method with an empirical application.
Section \ref{sec:con} concludes.
All technical proofs and auxiliary results are provided in the Supplementary Material.

\section{Identification of Area-specific Average Treatment Effect}\label{sec:set}

%

\subsection{Settings and estimand}

We consider a finite population consisting of $N$ individuals indexed by $i=1,\dots,N$. 
Each individual belongs to one of $J$ geographic areas, and the area to which individual $i$ belongs is denoted by $A_i \in \{1,\dots,J\}$. 
This notation differs slightly from the conventional notation in small area estimation, where units are often indexed separately within each area. 
However, the data structure considered here is essentially the same: each population unit belongs to exactly one area, and area-specific quantities are obtained by conditioning on, or aggregating over, units with the corresponding area label. 
The use of the common population index $i=1,\ldots,N$ and the area indicator $A_i$ is mainly for notational convenience, as it allows us to describe sampling, treatment assignment, and area membership within a unified individual-level framework. 
In particular, the target area-specific estimand is defined by conditioning on $A_i=j$, rather than by introducing a separate index system for each area.

For each individual we define a binary treatment indicator $T_i \in \{0,1\}$ and potential outcomes $Y_i(1)$ and $Y_i(0)$ under treatment and control, respectively. 
The observed outcome is
$$
Y_i = T_i Y_i(1) + (1-T_i)Y_i(0),
$$
following the potential outcomes framework \citep{R74, RR83}.

Our goal is to estimate area-specific average treatment effect (ATE) using survey data drawn from this population. 
The area-specific ATE is defined as 
\begin{align}\label{eqn:ate}
\tau_j\equiv
\mathbb{\mathbb{E}}\bigl[Y_i(1)-Y_i(0) \mid A_i=j\bigr],
\qquad j=1,\ldots,J .
\end{align}
This estimand captures heterogeneity in treatment effects across areas, and corresponds to a conditional average treatment effect widely studied in the causal inference literature \citep{H98, HI03}.
The main challenge arises because outcome, treatment status and several important covariates are observed only for survey respondents, while the target estimand is defined for the entire population of the target area $j$.

Regarding covariates, let $\X_i$ denote a vector of demographic covariates that are typically available for the population or from auxiliary sources such as census data.
Examples include age, sex, race, education, and income.
In addition, let $\Z_i$ denote individual-level characteristics that are observed only in the survey sample.
These variables capture attributes that cannot usually be obtained from administrative or census data.
In the empirical analysis in Section \ref{sec:rda}, they include political predispositions such as party identification, ideological self-placement, political interest.
Distinguishing between $\X_i$ and $\Z_i$ is important because $\Z_i$, although available only for sampled units, contain rich individual-level information that enables identification and improves the efficiency of population-level causal effect estimation.

\subsection{Sampled data}

Let $S_i \in \{0,1\}$ denote a sampling indicator such that $S_i=1$ if individual $i$ is included in the survey sample and
$S_i=0$ otherwise.
The observed data therefore consist of $\{Y_i,T_i,\X_i,\Z_i,A_i\}_{i=1; S_i=1}^n$.
Let $\pi_S(\X_i,A_i) = \mathbb{\mathbb{P}}(S_i=1| \X_i,A_i)$ denote the sampling probability.
This probability may vary across areas and demographic groups.
The sampling mechanism plays an important role because the survey sample is generally not a simple random sample of the population, as discussed in the survey sampling literature on informative sampling \citep{B83, PS99}.
Consequently, recovering area-specific causal effects requires accounting for the sampling process.


We consider several probabilities related to the distribution of areas in the population and in the survey sample.
First, define the marginal area probability $p_i(j) = \mathbb{\mathbb{P}}(A_i=j)$.
This quantity represents the population share of area $j$.
In many applications it can be obtained from external sources such as census data.
Second, define the area assignment probability conditional on demographic covariates $\pi_A(j | \X_i)=\mathbb{\mathbb{P}}(A_i=j | \X_i)$.
This probability describes how individuals with demographic characteristics $\X_i$ are distributed across areas.
It plays a role similar to the poststratification distribution in multilevel regression and poststratification, a widely used approach for combining survey and population information \citep{GL97}, to recover area-specific population quantities from survey data.
Finally, we also consider the conditional distribution of the area indicator within the survey sample, $\pi_A(j | \X_i,\Z_i)=\mathbb{\mathbb{P}}(A_i=j | \X_i,\Z_i,S_i=1)$.
This quantity describes the distribution of areas among sampled individuals conditional on their observed covariates.
It will appear in the estimation procedure when modeling the area assignment mechanism using survey data.
These three probabilities characterize different aspects of the area distribution: the marginal population distribution $\mathbb{\mathbb{P}}(A=j)$, the population assignment mechanism $\mathbb{\mathbb{P}}(A=j| \X)$, and the conditional distribution within the survey sample $\mathbb{\mathbb{P}}(A=j| \X,\Z,S=1)$.
Distinguishing these quantities is important for constructing estimators that combine survey data with auxiliary population information.
In practice, the population distribution $\mathbb{\mathbb{P}}(A=j|\X)$ may be obtained from auxiliary population data, whereas $\mathbb{\mathbb{P}}(A=j|\X,\Z,S=1)$ can be estimated from the survey sample.
Treatment is assumed to occur after the realization of the covariates and the sampling process.
Following \citet{RR83}, define the propensity score $e(\X_i,\Z_i)=\mathbb{\mathbb{P}}(T_i=1| \X_i,\Z_i,S_i=1)$.
This quantity represents the conditional probability of treatment assignment given the observed covariates among sampled units, and plays a central role in adjusting for selection into treatment in the observed sample.

Under the above setup, the full population consists of $\{Y_i(1),Y_i(0),T_i,\X_i,\Z_i,A_i\}_{i=1}^N$, but only the subset $\{Y_i,T_i,\X_i,\Z_i,A_i\}_{i=1; S_i=1}^n$ is observed.
Our objective is to recover the area-specific causal effects $\tau_j$ using the survey data together with auxiliary information about the population distribution of $\X$ and the area assignment mechanism.
A key feature of our setting is that the survey collects individual-specific covariates $\Z_i$ that are not available in population data sources such as the census, a setting related to problems with partially observed confounders in causal inference \citep{LR19}.
These variables may include characteristics such as political attitudes, policy preferences, or other behavioral attributes
that are important confounders of treatment assignment.
Consequently, valid adjustment for treatment assignment requires conditioning on both $\X_i$ and $\Z_i$.
While $\Z_i$ is observed only for sampled individuals ($S_i=1$), information about the population distribution of $\X_i$ and the area assignment mechanism allows us to combine survey and auxiliary data to recover the desired area-specific causal effects.
The framework developed in this paper shows how such survey-only covariates can be incorporated into causal small area estimation.

In the next section, we present identification results that express $\tau_j$ in terms of observable quantities under suitable assumptions.
This setting arises naturally in many applications where individual attitudes or behavioral characteristics are measured
only in surveys but play an important role in treatment assignment.

\subsection{Identification: assumptions and results}
\label{sec:iden}

We begin with the conventional direct identification formula based on within-area observations, which can be unstable in small areas. 
We then develop an alternative identification representation that incorporates information across areas while preserving the target parameter. 
We introduce the assumptions required for identifying the area-specific causal effects defined in \eqref{eqn:ate}.

\begin{itemize}
\item[(A1)]
(Sampling ignorability) \ 
$(Y_i(1),Y_i(0),T_i,\Z_i)\perp S_i\mid (\X_i,A_i)$.

\item[(A2)]
(Unconfoundedness) \ $(Y_i(1),Y_i(0))\perp T_i\mid (\X_i,\Z_i,S_i=1)$.

\item[(A3)]
(Overlap) \ There exists a constant $\epsilon>0$ such that three probabilities, $\mathbb{P}(S_i=1| \X_i,A_i=j)$, $\mathbb{P}(A_i=j| \X_i,\Z_i,S=1)$ and $\mathbb{P}(T_i=1| \X_i,\Z_i,S_i=1)$ lie in $(\epsilon, 1-\epsilon)$ for all $j$ and all values of the covariates.

\item[(A4)]
(Area ignorability) \ $(Y_i(1),Y_i(0),T_i)\perp A_i\mid(\X_i,\Z_i,S_i=1)$.
\end{itemize}

Assumptions~(A1)-(A3) are relatively standard ones. 
Assumption~(A1) states that conditional on demographic covariates $\X_i$ and the area indicator $A_i$, inclusion in the survey sample is independent of the potential outcomes, treatment status, and survey-specific covariates.
It implies that the survey sample is representative of the population within strata defined by $(\X,A)$, which is standard in analyses of survey data under ignorable sampling designs \citep{S03, B83}.
Assumption (A2) states that treatment assignment is independent of the potential outcomes conditional on the observed covariates $(\X,\Z)$ among sampled individuals, and is standard in observational causal inference \citep{RR83, IR15}.
Assumption (A3) ensures that each sampling status, each area and each treatment level occurs with positive probability within the relevant covariate strata, which is a standard in causal inference \citep{IR15}.

On the other hand, Assumption~(A4) is less standard and is analogous to the area ignorability condition considered in \citet{KY21}.
It states that, among sampled individuals, the potential outcomes and treatment status are independent of the area indicator after conditioning on the observed covariates $(\X_i,\Z_i)$.
Thus, individuals with the same demographic characteristics $\X_i$ and survey-only attributes $\Z_i$ are assumed to have the same joint distribution of potential outcomes and treatment status, regardless of their area.
The survey-only covariates $\Z_i$ are crucial in this assumption because they may capture individual-level characteristics, such as political attitudes or behavioral traits, that explain systematic differences across areas beyond those captured by $\X_i$ alone.
This interpretation is conceptually related to M-quantile SAE \citep[e.g.][]{chambers2006m,schirripa2026bias}, which represents between-area heterogeneity through covariate-adjusted features of the conditional outcome distribution rather than relying solely on explicit area-specific random effects.
Similarly, Assumption~(A4) reflects the idea that area differences relevant to potential outcomes and treatment assignment can be sufficiently explained by rich observed covariates.
In our framework, the area assignment probability $\pi_A(j\mid \X_i,\Z_i)$ summarizes how these covariate profiles are distributed across areas and provides the basis for transporting information from the full survey sample to the target area.

Under these assumptions, the area-specific treatment effects introduced in \eqref{eqn:ate} can be expressed in terms of observable quantities.
A natural identification strategy is to rely solely on observations within area $j$. 
Under Assumptions (A1)-(A3), by using inverse probability weighting arguments standard in causal inference \citep{HT52, RR83}, one can express $\tau_j$ using only units with $A=j$ as
\begin{align}\label{eqn:de}
\tau_j =
\mathbb{\mathbb{E}}\left[
\frac{\mathbf{1}\{S_i=1,A_i=j\}}
{\pi_S(\X_i,A_i=j) p_i(j)}
\left\{
\frac{T_iY_i}{e_j(\X_i,\Z_i)}
-
\frac{(1-T_i)Y_i}{1-e_j(\X_i,\Z_i)}
\right\}
\right],
\end{align}
where $e_j(\X_i,\Z_i)=\mathbb{P}(T_i=1| \X_i,\Z_i,A_i=j,S_i=1)$, $\pi_S(\X_i,A_i)=\mathbb{P}(S_i=1| \X_i,A_i)$ and $p_i(j)=\mathbb{P}(A_i=j)$.

While this representation is conceptually simple, it is unsuitable for small area settings. 
When the number of sampled units in each area is limited, direct estimators based on this expression exhibit high variance and can be unstable in practice. 
To address the instability of the direct estimator, we introduce an identification formula that incorporates observations from all areas in the following proposition.

\begin{prp}\label{prp:id}
Under Assumptions (A1)-(A4), the area-specific average treatment effect, $\tau_j$ can be expressed as 
\begin{align*}
\tau_j=&
\mathbb{\mathbb{E}}\left[
\frac{\one\{S_i=1\}}{\pi_S(\X_i,A_i=j)} \frac{\pi_A(j |\X_i, \Z_i)}{p_i(j)}
\left\{
\frac{T_iY_i}{e(\X_i,\Z_i)}
-
\frac{(1-T_i)Y_i}{1-e(\X_i,\Z_i)}
\right\}
\right]
\end{align*}
where $e(\X_i,\Z_i)=\mathbb{P}(T_i=1| \X_i,\Z_i,S_i=1)$, $\pi_A(j | \X_i,\Z_i)=\mathbb{P}(A_i=j | \X_i,\Z_i,S_i=1)$, $\pi_S(\X_i,A_i)=\mathbb{P}(S_i=1| \X_i,A_i)$, and $p_i(j)=\mathbb{P}(A_i=j)$ for $j=1,\ldots,J$.
\end{prp}

This type of reweighting representation is closely related to transportability and covariate shift adjustments studied in the causal inference literature \citep{DH19, T20}.
This representation reweights observations from all areas to recover the distribution of $(\X,\Z)$ in area $j$. 
The factor $\pi_A(j| \X,\Z)$ enables cross-area transport, while $\pi_S(\X,A=j)$ adjusts for the sampling mechanism. 
Although this representation improves stability relative to the direct estimator, it relies entirely on correct specification of the weighting components and is therefore not robust.

The next corollary shows that the identification formula in Proposition \ref{prp:id} admits a decomposition that reveals how information is combined across areas.

\begin{cor}\label{cor:aif}
The area-specific average treatment effect is also identified as
\begin{align*}
\tau_j&=
\mathbb{\mathbb{E}}\left[\frac{\one\{S_i=1,A_i=j\}}{\pi_S(\X_i,A_i=j)p_i(j)}
\left\{
\frac{T_iY_i}{e(\X_i,\Z_i)}
-
\frac{(1-T_i)Y_i}{1-e(\X_i,\Z_i)}
\right\}
\pi_A(j | \X_i,\Z_i) \right] \\
&+\mathbb{\mathbb{E}}\left[\frac{\one\{S_i=1,A_i\neq j\}}{\pi_S(\X_i,A_i=j)p_i(j)}
\left\{
\frac{T_iY_i}{e(\X_i,\Z_i)}
-
\frac{(1-T_i)Y_i}{1-e(\X_i,\Z_i)}
\right\}  \pi_A(j | \X_i,\Z_i)
\right].
\end{align*}
\end{cor}

The first term corresponds closely to the direct estimator, which relies solely on observations from the target area, which is typically unstable in small samples. 
The second term represents a partial pooling component that incorporates information from other areas through the weighting structure using the area assignment mechanism $\pi_A(j| \X,\Z)$.
This component plays a crucial role when the number of sampled individuals in the target area is limited, allowing the estimator to borrow strength from observations in other areas while adjusting for differences in covariate distributions.

The structure of this decomposition closely resembles the identification framework developed by \citet{KY21} for small area mean estimation under sampling and area ignorability.
In particular, both approaches combine a within-area contribution with an additional term that integrates information from other areas through covariate adjustment.
Our result extends this idea to the causal inference setting by incorporating inverse probability weighting for treatment
assignment.

\subsection{Doubly robust identification}

We now present our main result, which combines outcome regression with the weighting representation to obtain a doubly robust characterization.
Define outcome regression model $m_t(\X,\Z)$ as
$$
m_t(\X,\Z) = \mathbb{\mathbb{E}}[Y | \X,\Z,T=t,S=1], \quad t=0,1.
$$
Let $\bta$ be a set of nuisance functions of outcome regression models, propensity score and area assignment probability, namely $\bta=(m_1,m_0,e,\pi_A)$.
For $\W_i=(Y_i,\X_i,\Z_i,T_i,A_i,S_i)$, define
\begin{align}\label{eqn:if}
\begin{split}
&\phi_{ij}(\W_i, \bta) 
= 
\frac{\one\{S_i=1\}}{\pi_S(\X_i,A_i=j)} \frac{\pi_A(j|\X_i,\Z_i)}{p_i(j)} \\
& \times \left[ \frac{T_i\{Y_i-m_1(\X_i,\Z_i)\}}{e(\X_i,\Z_i)}-
\frac{(1-T_i)\{Y_i-m_0(\X_i,\Z_i)\}}{1-e(\X_i,\Z_i)} 
+m_1(\X_i,\Z_i)-m_0(\X_i,\Z_i)
\right] \\
& +\frac{\one\{S_i=1\}}{\pi_S(\X_i,A_i=j)}\frac{\one\{A_i=j\}-\pi_A(j|\X_i,\Z_i)}{p_i(j)} \{m_1(\X_i,\Z_i)-m_0(\X_i,\Z_i)\}.
\end{split}
\end{align}
Let $\widetilde{\bta}=(\widetilde{m}_1,\widetilde{m}_0,\widetilde{e},\widetilde{\pi}_A)$ be arbitrary models for the true, unknown nuisance functions $\bta=(m_1,m_0,e,\pi_A)$.
The next theorem shows that the expectation of the estimating function in \eqref{eqn:if} satisfies a doubly robust property.
Such doubly robust representations are well known in causal inference and missing data problems \citep{RRZ94, BR05, T06}.

\begin{thm}\label{thm:dr}
Under Assumptions (A1)-(A4), $\tau_j = {\mathbb{E}}[\phi_{ij}(\W_i, \widetilde{\bta})]$
holds if either $\widetilde{m}_t = m_t$ for $t=0,1$ or $(\widetilde{e}, \widetilde{\pi}_A) = (e, \pi_A)$ is satisfied.
\end{thm}

The above representation integrates three sources of selection: treatment assignment, survey sampling, and area assignment, extending standard doubly robust estimators to settings with additional selection mechanisms \citep{BR05}.
The second term in the estimating function in \eqref{eqn:if} corrects for discrepancies in the area distribution and has no analogue in standard doubly robust estimators. 
This term is essential for leveraging information across areas while preserving identification of $\tau_j$.

As shown in the proof of Theorem \ref{thm:dr} in the Supplementary Material, in the case where the outcome regression model is correctly specified but the area assignment model is misspecified, which satisfies the doubly robust condition in Theorem \ref{thm:dr}, the estimand reduces to the following equation
\begin{align*}
\tau_j
=
\mathbb{\mathbb{E}}\bigg[\frac{\one\{S_i=1, A_i=j\}}{\pi_S(\X_i,A_i=j)p_i(j)} 
\{m_1(\X_i,\Z_i) -m_0(\X_i,\Z_i)\} \bigg].
\end{align*}
In contrast to Proposition \ref{prp:id}, the above expression depends only on individuals in area $j$.
However, as discussed in Section \ref{sec:est}, the correctly specified outcome regression models are estimated using observations from all areas, thereby enabling borrowing strength across areas and leading to more stable estimation of area-specific ATE than the doubly robust estimand based on the direct estimate in \eqref{eqn:de}.

We now establish that the proposed score corresponds to the efficient influence function in a nonparametric model, following the general theory of semiparametric efficiency \citep{BR05, V00}.
Next theorem shows that the score given in ($\ref{eqn:if}$) can achieve the semiparametric efficiency bound for the ATE under Assumptions (A1)–(A4).

\begin{thm}\label{thm:seb}
Suppose Assumptions (A1)-(A4) hold.
The efficient influence function for the area-specific ATE, $\tau_j$ is given by $\varphi_{ij}(\W_i, \bta_0)= \phi_{ij}(\W_i,\bta_0)-\tau_j$, where $\bta_0=(m_{1,0},m_{0,0},e_0,\pi_{A,0})$ denotes the true nuisance functions, and the semiparametric efficiency bound for any regular estimators for $\tau_j$ is
\begin{align*}
&{\mathbb{E}}[\{\varphi_{ij}(\W_i, \bta_0)\}^2]  \\
&=
\mathbb{E}\Bigg[
\frac{\pi_{A,0}(j|\X_i,\Z_i)^2}{\pi_S(\X_i,A_i=j)p_i(j)^2} 
\Bigg\{\frac{\mathbb{V}(Y_i(1)|\X_i,\Z_i)}{e_0(\X_i,\Z_i)}+
\frac{\mathbb{V}(Y_i(0)|\X_i,\Z_i)}{1-e_0(\X_i,\Z_i)}\Bigg\}\Bigg] \\
& \quad +
\mathbb{E}\Bigg[
\frac{\pi_{A,0}(j|\X_i,\Z_i)}{\pi_S(\X_i,A_i=j)p_i(j)^2}
\{m_{1,0}(\X_i,\Z_i)-m_{0,0}(\X_i,\Z_i)-\tau_j\}^2\Bigg],
\end{align*}
where $\mathbb{V}(Y_i(t)|\X_i,\Z_i)=\mathbb{E}[\{Y_i(t)-m_{t,0}(\X_i,\Z_i)^2\}|\X_i,\Z_i,S_i=1]$ for $t=0,1$.
\end{thm}

The representation of the semiparametric efficiency bound in Theorem \ref{thm:seb} provides insight into the sources of estimation uncertainty.
In particular, the bound decomposes into two components.
The first captures outcome noise through the conditional variance of $Y_i$ given $(\X_i,\Z_i,T_i,S_i=1)$, scaled by the weighting factor $\pi_{A,0}(j|\X_i,\Z_i)^2 / \{\pi_S(\X_i,A_i=j)p_i(j)^2\}$, reflecting variability induced by inverse probability weighting.
The second captures treatment effect heterogeneity, measured by the dispersion of the conditional average treatment effect $m_{1,0}(\X_i,\Z_i)-m_{0,0}(\X_i,\Z_i)$ around $\tau_j$, representing intrinsic variation across individuals.

\section{Causal SAE: Estimation and Inference}\label{sec:est}

This section describes the estimation procedure for the proposed causal small area estimator.
The representation in \eqref{eqn:if} naturally suggests a plug-in estimation strategy based on estimated nuisance functions.
This plug-in strategy is standard in semiparametric causal inference and is closely related to doubly robust estimation methods \citep{RRZ94, BR05, T06}.
Specifically, estimation requires three components: the outcome regression functions, the propensity score, and the area assignment probability. 
Once these quantities are estimated, they can be substituted into the sample analogue of the doubly robust score function to obtain estimators of the area-specific average treatment effects.

\subsection{Estimation of nuisance functions}

The following nuisance functions are required to be estimated for the proposed estimators:
\begin{align*}
&(\text{Outcome regression}) \quad  m_t(\X,\Z)\equiv \mathbb{\mathbb{E}}[Y|\X,\Z,T=t,S=1], \quad t=0,1,\\
&(\text{Propensity score}) \quad e(\X,\Z)\equiv \mathbb{P}(T=1|\X,\Z,S=1)\\
&(\text{Area assignment probability}) \quad \pi_A(j | \X,\Z)\equiv \mathbb{P}(A=j | \X,\Z,S=1)
\end{align*}
Although the framework allows flexible semiparametric or nonparametric estimators, we estimate these quantities using parametric working models for simplicity and transparency.
In particular, we consider the following parametric specifications: 
\begin{align*}
&m_t(\X,\Z;\bbe_t)=\X' \bbe_{t1} + \Z' \bbe_{t2}, \quad e(\X,\Z;\bga)
=\frac{\exp(\X' \bga_1 + \Z' \bga_2)}
{1+\exp(\X' \bga_1 + \Z' \bga_2)},\\ 
&\pi_A(j|\X,\Z;\bal)=
\frac{\exp(\X' \bal_{1,j}+\Z' \bal_{2,j})}
{\sum_{\ell=1}^{m}\exp(\X' \bal_{1,\ell}+\Z' \bal_{2,\ell})},
\end{align*}
where $\bbe_t=(\bbe_{t1}',\bbe_{t2}')'$ for $t=0,1$, $\bga=(\bga_1',\bga_2')'$ and $\bal=(\bal_{1,1}',\bal_{2,1}',\ldots,\bal_{1,J}',\bal_{2,J}')'$ are unknown parameter vectors.
Then, we estimate $\bbe_t$ by ordinary least squares using survey observations with $T_i=t\in\{0,1\}$.
The parameters $\bga$ and $\bal$ are estimated by maximum likelihood using all survey observations, $\{Y_i, T_i, \X_i, \Z_i, A_i\}_{i=1; S_i=1}^n$.
Such parametric specifications are commonly used in applied causal inference and survey analysis \citep{MN89, M86}.

Moreover, to mitigate overfitting bias and ensure valid asymptotic inference of treatment effect, we employ cross-fitting. 
This approach follows recent developments in double/debiased machine learning \citep{CCD18}.
Specifically, we partition the sample into $K$ mutually exclusive folds $\{I_k\}_{k=1}^K$ of approximately equal size.
For each fold $k$, we estimate the nuisance parameters using observations not in $I_k$, denoted by $I_k^c$ and let $\widehat{\bbe}_t^{(-k)}$, $\widehat{\bga}^{(-k)}$ and $\widehat{\bal}^{(-k)}$ denote the resulting estimators.
We then evaluate the nuisance functions on observations in $I_k$ as predictions obtained from models trained on $I_k^c$.
Namely, for each observation $i \in I_k$, the estimated nuisance functions are $m_t(\X_i,\Z_i; \bbeh_t^{(-k)})$, $e(\X_i,\Z_i; \bgah^{(-k)})$ and $\pi_A(j | \X_i,\Z_i; \balh^{(-k)})$.
Then, by combining across folds, the cross-fitted nuisance estimators are defined as
\begin{align}\label{eqn:ne}
&\widehat{m}_t(\X_i,\Z_i)
=
\sum_{k=1}^K
\mathbf{1}\{i \in I_k\}
m_t(\X_i,\Z_i;\bbeh_t^{(-k)}), \nonumber\\
&\widehat{e}(\X_i,\Z_i)
=
\sum_{k=1}^K
\mathbf{1}\{i \in I_k\}
e(\X_i,\Z_i; \bgah^{(-k)}), \\
&\widehat{\pi}_A(j| \X_i,\Z_i)
=
\sum_{k=1}^K
\mathbf{1}\{i \in I_k\}
\pi_A(j | \X_i,\Z_i; \balh^{(-k)}). \nonumber
\end{align}

Cross-fitting ensures that the nuisance estimates used for each observation are obtained from independent data. 
This property plays a crucial role in eliminating the first-order bias induced by nuisance estimation errors, which allows for valid asymptotic inference even if flexible machine learning methods are used.
In particular, as shown in the next subsection, the resulting estimators admit a simple asymptotic linear representation without requiring additional correction terms for nuisance estimation.

\subsection{Auxiliary probabilities}

The final estimator given in Corollary~\ref{cor:aif} involve a probability $\pi_S(\X,A=j)$ that depend on the sampling design.
We first note that 
$$
\pi_S(\X,A=j)\equiv \mathbb{P}(S=1|\X,A=j)=\frac{\mathbb{P}(S=1|\X)\mathbb{P}(A=j|\X,S=1)}{\mathbb{P}(A=j|\X)}.
$$
Then, the sampling probability $\mathbb{P}(S=1|\X)$ can be evaluated using survey weights.
In practice, this quantity can be approximated using the inverse of survey weights, for example, calibration or raking weights constructed to align the sample with known population margins.
See \citet{DS92} and \citet{S03} for discussions of calibration and weighting methods in survey sampling.
The probability $\mathbb{P}(A=j|\X,S=1)$ can be estimated nonparametrically from the survey data, while $\mathbb{P}(A=j|\X)$ can be evaluated using auxiliary population data such as census information.

\subsection{Estimation of the Average Treatment Effect}
\label{subsec:est}

We propose estimators for the area-specific average treatment effect. 
Throughout, we use the cross-fitted nuisance estimators introduced in Section 4.1.
We rewrite the score function in \eqref{eqn:if} coinciding with the efficient influence function for $\tau_j$ as
$$
\phi_{ij}=
\phi_{ij}(\W_i;\bta)=
w_{ij1} \phi_{i1} + w_{ij2} \phi_{i2},
$$
where
\begin{align*}
&\phi_{i1}
=
\frac{T_i\{Y_i-m_1(\X_i,\Z_i;\bbe_1)\}}{e(\X_i,\Z_i;\bga)}
-
\frac{(1-T_i)\{Y_i-m_0(\X_i,\Z_i;\bbe_0)\}}{1-e(\X_i,\Z_i;\bga)} \\
&\qquad +
m_1(\X_i,\Z_i;\bbe_1)-m_0(\X_i,\Z_i;\bbe_0), \\
&\phi_{i2}
=
m_1(\X_i,\Z_i;\bbe_1)-m_0(\X_i,\Z_i;\bbe_0),
\end{align*}
and
$$
w_{ij1}
=
\frac{\one\{S_i=1\}}{p_i(j)}
\frac{\pi_A(j|\X_i,\Z_i;\bal)}{\pi_S(\X_i,A_i=j)},
\qquad
w_{ij2}
=
\frac{\one\{S_i=1\}}{p_i(j)}
\frac{\one\{A_i=j\}-\pi_A(j|\X_i,\Z_i;\bal)}{\pi_S(\X,A=j)}.
$$
Let $\widehat{\phi}_{ik}$ and $\widehat{w}_{ijk}$ for $k=1,2$ be the estimator of $\phi_{ik}$ and $w_{ijk}$ for $k=1,2$, respectively, obtained by replacing nuisance functions $\bta$ with their estimators $\widehat{\bta}=(\widehat{m}_1, \widehat{m}_0, \widehat{e}, \widehat{\pi}_A)$, and let 
\begin{align}\label{eqn:phih}
\widehat{\phi}_{ij}=\phi_{ij}(\W_i;\widehat{\bta})=\widehat{w}_{ij1} \widehat{\phi}_{i1} + \widehat{w}_{ij2} \widehat{\phi}_{i2}.
\end{align}
Then, the following two estimators are extensions of classical Horvitz--Thompson and H\'ajek estimators \citep{HT52, H71}, respectively., defined as 
\begin{align}\label{eqn:HT}
\widehat{\tau}^{\rm HT}_j
=
\frac{1}{n} \sum_{i=1}^n \widehat{\phi}_{ij}, \qquad 
\widehat{\tau}^{\rm Hajek}_j
=
\frac{\sum_{i=1}^n \widehat{\phi}_{ij}}{\sum_{i=1}^n \widehat{w}_{ij1}}.
\end{align}
The use of cross-fitting, together with the orthogonal structure of the score, yields an estimator that admits an asymptotically linear representation.
This representation forms the basis for the asymptotic theory developed in the subsequent section.

The overall estimation procedure can be summarized as follows.

\begin{algo}[Estimation of area-specific ATEs]
\end{algo}

\begin{enumerate}
\item 
Split the survey sample ($S=1$) into $K$ mutually exclusive folds $\{\mathcal{I}_k\}_{k=1}^K$.
For each fold $k$, estimate nuisance functions using observations not in $\mathcal{I}_k$:

\begin{itemize}
\item 
For outcome regression model, $\bbeh_t^{(-k)}$ is estimated by ordinary least squares.

\item 
For propensity score, $\bgah^{(-k)}$ is estimated by maximum likelihood.

\item 
For area assignment probability, $\balh^{(-k)}$ is estimated by maximum likelihood.
\end{itemize}

Define $\widehat{\bta}=(\widehat{m}_0, \widehat{m}_1, \widehat{e}, \widehat{\pi}_A)$ given in \eqref{eqn:ne}.

\item 
Compute sampling probability $\pi_S(\X,A=j)=\mathbb{P}(S=1|\X,A=j)$ by obtaining estimates of the auxiliary probabilities $\mathbb{P}(S=1|\X,A)$, $\mathbb{P}(A=j|\X,S=1)$, and $\mathbb{P}(A=j|\X)$ using survey weights and population data.

\item 
Compute $\widehat{\phi}_{ij}=\phi_{ij}(\W_i;\widehat{\bta})$, given in \eqref{eqn:phih}.

\item 
Compute the Horvitz--Thompson type estimator $\hat{\tau}^{\rm HT}_j$ and the H\'ajek type estimator $\hat{\tau}^{\rm Hajek}_j$ defined in \eqref{eqn:HT}.

\end{enumerate}

\subsection{Asymptotic properties and variance estimation}

In this subsection, we derive feasible estimators for the asymptotic variance of the proposed ATE estimators based on the influence function representation established in Section \ref{subsec:est}.
Let $\bta_0$ denote the true nuisance functions $\bta$, and $\mu_{wj1} = \mathbb{\mathbb{E}}[w_{ij1}]$.
Under additional regularity conditions in the Supplementary Material and cross-fitting, the ATE estimators given in (\ref{eqn:HT}) have an asymptotic linear representation given in the following theorem.
Such representations are central in semiparametric efficiency theory \citep{BR93, V00}.

\begin{thm}\label{thm:an}
Under Assumptions (A1)-(A4) in Section~\ref{sec:iden} and Assumptions (SA1)-(SA3) in the Supplementary Material, the Horvitz--Thompson type estimator and the H\'ajek type estimator have the following asymptotic linear representations under $n\to\infty$: 
\begin{align*}
&\sqrt{n}(\widehat{\tau}^{\rm HT}_j - \tau_j)
= \frac{1}{\sqrt{n}} \sum_{i=1}^n \{ \phi_{ij}(\W_i;\bta_0) - \tau_j \} + o_p(1),\\
&\sqrt{n}(\widehat{\tau}^{\rm Hajek}_j - \tau_j)
= \frac{1}{\sqrt{n}} \sum_{i=1}^n \frac{\phi_{ij}(\W_i;\bta_0) - \tau_j w_{ij1}}{\mu_{wj1}}
+ o_p(1).
\end{align*}

\end{thm}

\medskip
From Theorem \ref{thm:an}, it is clear that the asymptotic variances of $\widehat{\tau}^{\rm HT}$ and $\widehat{\tau}^{\rm Hajek}$ are respectively given by
\begin{align*}
&\Var(\widehat{\tau}^{\rm HT})=
\frac{1}{n}
\mathbb{\mathbb{E}}
\left[
\left\{
\phi_{ij}(\W_i; \bta_0) - \tau_j
\right\}^2
\right]
+o(n^{-1}),\\
&\Var(\widehat{\tau}^{\rm Hajek})=
\frac{1}{n\mu_{wj1}^2}
\mathbb{\mathbb{E}}
\left[
\left\{ \phi_{ij}(\W_i; \bta_0) - \tau_j w_{ij1}\right\}^2 
\right]
+o(n^{-1}).
\end{align*}
As a consequence of the influence function representation, these asymptotic variances coincide with the semiparametric efficiency bounds, which the proposed estimators of area-specific ATE in \eqref{eqn:HT} can achieve when all nuisance functions are consistently estimated.

Based on the above influence function representations, we estimate the variances by their sample analogues, defined as 
\begin{align}
\label{eqn:vht}
&\widehat{V}^{\rm HT}_j= \frac{1}{n^2} \sum_{i=1}^n (\widehat{\phi}_{ij} - \widehat{\tau}^{\rm HT}_j)^2, \qquad 
\widehat{V}_j^{\rm Hajek}=
\frac{1}{n^2\overline{w}_{j1}^2}
\sum_{i=1}^n
\left(\widehat{\phi}_{ij} - \widehat{\tau}_j^{\rm Hajek} \widehat{w}_{ij1}\right)^2,
\end{align}
where $\overline{w}_{j1} = n^{-1} \sum_{i=1}^n \widehat{w}_{ij1}$.
These variance estimators remain valid despite the use of estimated nuisance functions. 
The key reason is that cross-fitting ensures that, for each observation, the nuisance estimators are obtained from an independent subsample.
As a result, the first-order expansion is unaffected by the estimation error in nuisance parameters, and its contribution is asymptotically negligible, namely $n^{-1/2} \sum_{i=1}^n \{ \widehat{\phi}_{ij} - \phi_{ij}(\X_i,\Z_i;\bta_0) \}= o_p(1)$.
This greatly simplifies inference compared to the standard plug-in approach without sample splitting.

\section{Simulation Study}\label{sec:sim}

We conduct a set of Monte Carlo experiments to evaluate the finite-sample performance of the proposed estimators.
The simulation study is designed to evaluate three aspects of the proposed estimators: 
(i) stability relative to the direct estimator, 
(ii) efficiency gains implied by the efficient influence function representation, and 
(iii) robustness to misspecification of nuisance functions.

\subsection{Settings and estimators}
Consider a population consisting of $J=50$ areas with two scenarios of total population size, $N=400{,}000$ and $N=800{,}000$.
For each individual $i(=1,\ldots,N)$, we generate covariates $X_i, Z_i$ by $X_i \sim N(0,1)$ and $Z_i \sim N(\mu_{z,i},1)$ for $\mu_{z,i}\sim DU(1,50)/10$, where $DU(a,b)$ denotes discrete uniform distribution on $[a,b]$.
For generic covariates $\D=(D_1,D_2)$, outcome regression models, propensity score and area assignment probability are respectively set as follows; 
\begin{align*}
&f_0(\D)=1+0.5D_1+0.7D_2, \quad f_1(\D)=2+D_1+1.4D_2, \\
&f_T(\D) = \exp(-0.2+0.4D_1+0.4D_2)/\{1+\exp(-0.2+0.4D_1+0.4D_2)\}, \\
&f_A(\D) = \exp(\al_{j1}+\al_{j2}D_1+\al_{j2}D_2)/\sum_{j'=1}^m\exp(\al_{j'1}+\al_{j'2}D_1+\al_{j'2}D_2)\},
\end{align*}
where it is assumed to be $T\perp S|\D$ and $A\perp S|\D$, and $\al_{j,k}\sim N(0,0.15^2)$ for $j=1,\ldots,50$ and $k=1,2,3$.

We consider eight types of data generating process (DGP), which is partly inspired by simulation settings commonly used in the literature on doubly robust estimation, including \citet{SZ20}.
Each DGP is defined by changing whether each component of the DGP depends on the original covariates $(X,Z)$ or the transformed covariates $(X^2,Z^2)$. 
Specifically, for each DGP, the potential outcomes, propensity score, and area assignment probability are generated as $Y(t)= f_t(\D_Y)+\ep_t$ with $t=0,1$, $\mathbb{P}(T=1|X,Z)= f_T(\D_T)$ and $\mathbb{P}(A=j|X,Z)= f_A(\D_A)$, where $\D_Y,\D_T,\D_A$ are specified in Table~\ref{tab:dgp}.
In all DGPs, the error terms are independently generated as $\ep_0,\ep_1\sim N(0,10^2)$.
From the population, a survey sample is drawn according to the sampling mechanism
$$
\mathbb{P}(S_i=1|X_i,A_i) = \mathbb{P}(S_i=1|X_i) = 1/\{1+\exp(4-0.3X_i)\},
$$
which is assumed to be known.
Let $Y_i=T_iY_i(1)+(1-T_i)Y_i(0)$.
Assume that $(X_i,A_i)$ are observed for all subjects and $(Y_i,Z_i,T_i)$ are observed only for subjects with $S_i=1$.
Then, DGP 1--5 satisfy the doubly robustness condition in Theorem \ref{thm:dr}.

\begin{table}[t]
\centering
\caption{Covariates used in each component of the data-generating process in the simulation study.}
\label{tab:dgp}
\begin{tabular}{c|ccc}
\hline
DGP & Outcome regression $\D_Y$ & Propensity score $\D_T$ & Area assignment $\D_A$ \\
\hline
1 & $(X,Z)$ & $(X,Z)$ & $(X,Z)$ \\
2 & $(X,Z)$ & $(X^2,Z^2)$ & $(X,Z)$ \\
3 & $(X,Z)$ & $(X,Z)$ & $(X^2,Z^2)$ \\
4 & $(X,Z)$ & $(X^2,Z^2)$ & $(X^2,Z^2)$ \\
5 & $(X^2,Z^2)$ & $(X,Z)$ & $(X,Z)$ \\
6 & $(X^2,Z^2)$ & $(X,Z)$ & $(X^2,Z^2)$ \\
7 & $(X^2,Z^2)$ & $(X^2,Z^2)$ & $(X,Z)$ \\
8 & $(X^2,Z^2)$ & $(X^2,Z^2)$ & $(X^2,Z^2)$ \\
\hline
\end{tabular}
\end{table}

We compare the Horvitz--Thompson type estimator $\widehat{\tau}^{\rm HT}_j$ and the H\'ajek type estimator $\widehat{\tau}^{\rm Hajek}_j$ defined in \eqref{eqn:HT}. 
Both estimators are constructed using estimated nuisance functions, including the outcome regression, the propensity score, and the area assignment probability. 
The nuisance functions are estimated using parametric working models, as described in Section \ref{sec:est}. 
For inference, the variances of $\widehat{\tau}^{\rm HT}_j$ and $\widehat{\tau}^{\rm Hajek}_j$ are estimated using the estimators given in \eqref{eqn:vht}, respectively.
As an additional benchmark, we consider a direct estimator based on the identification formula in \eqref{eqn:de}, modified to admit a doubly robust representation. 
Specifically, we construct its sample analogue using a H\'ajek type normalization, yielding a doubly robust direct estimator for each area, and we denote the estimator as $\widehat{\tau}_j^{Direct}$.
This benchmark allows us to assess the gains from incorporating information across areas relative to within-area estimation under a comparable doubly robust structure.

\subsection{Performance measures and results}
For each estimator of area-specific average treatment effect, we evaluate bias and root mean squared error as the performance criteria, which are respectively defined as
$$
\mathrm{Bias}(\widehat{\tau}_j) =
\frac{1}R \sum_{r=1}^{R} \widehat{\tau}_j^{(r)} - \tau_j,
\qquad
\mathrm{RMSE}(\widehat{\tau}_j) =
\sqrt{\frac{1}R \sum_{r=1}^{R} (\widehat{\tau}_j^{(r)} - \tau_j)^2},
$$
where $\widehat{\tau}_j^{(r)}$ denotes the estimate obtained in the $r$-th Monte Carlo replication. 
To quantify the efficiency gain, we can compute the percentage relative improvement in average loss (PRIAL) of the proposed estimator $\widehat{\tau}_j$ over the direct estimator $\widehat{\tau}_j^{\rm Direct}$, defined as
$$
\mathrm{PRIAL}(\widehat{\tau}_j) = 
100 \times \left\{1 - \frac{\mathrm{RMSE}(\widehat{\tau}_j)}{\mathrm{RMSE}(\widehat{\tau}_j^{\rm Direct})}\right\}.
$$
To evaluate the performance of variance estimation for the proposed method, we also consider the ratio of the estimated variance to the true variance, given by
$\mathrm{Ratio}(\widehat{V}_j)=\widehat{V}_j^{(r)}/\Var(\widehat{\tau}_j)$, where $\widehat{V}_j^{(r)}$ is the variance estimator associated with $\widehat{\tau}_j^{(r)}$, and the true variance is approximated by $\Var(\widehat{\tau}_j)=R^{-1} \sum_{r=1}^{R} \left(\widehat{\tau}_j^{(r)} - \tau_j\right)^2$.
Values close to one indicate accurate variance estimation.
All reported results are based on $R=2000$ Monte Carlo replications.

\begin{figure}[t]
  \centering
  \begin{minipage}{0.4\columnwidth}
     \centering
     \includegraphics[width=\columnwidth]{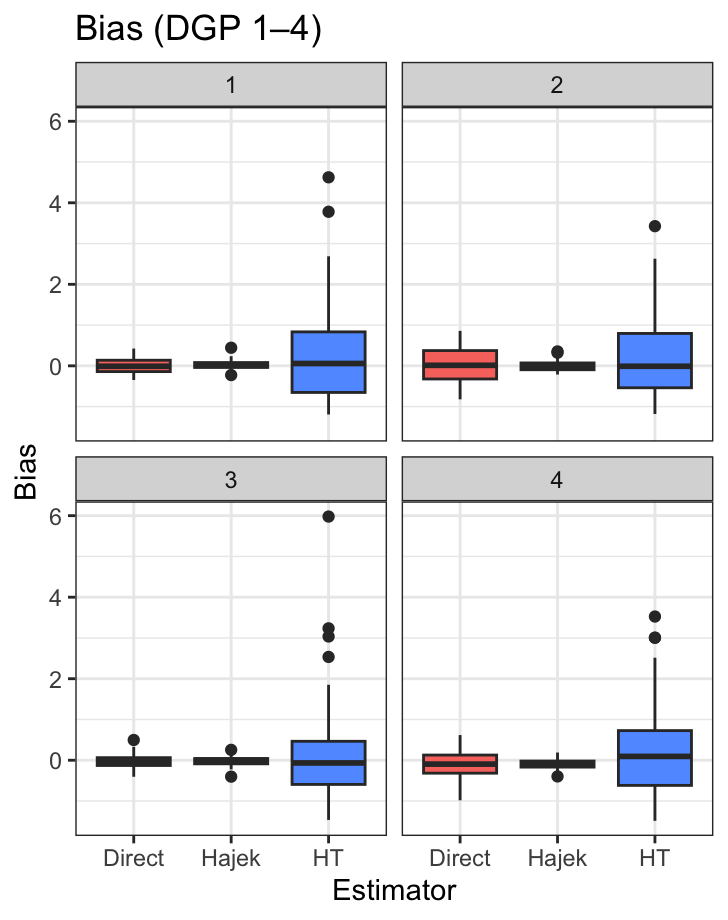}
  \end{minipage}
  \begin{minipage}{0.4\columnwidth}
     \centering
     \includegraphics[width=\columnwidth]{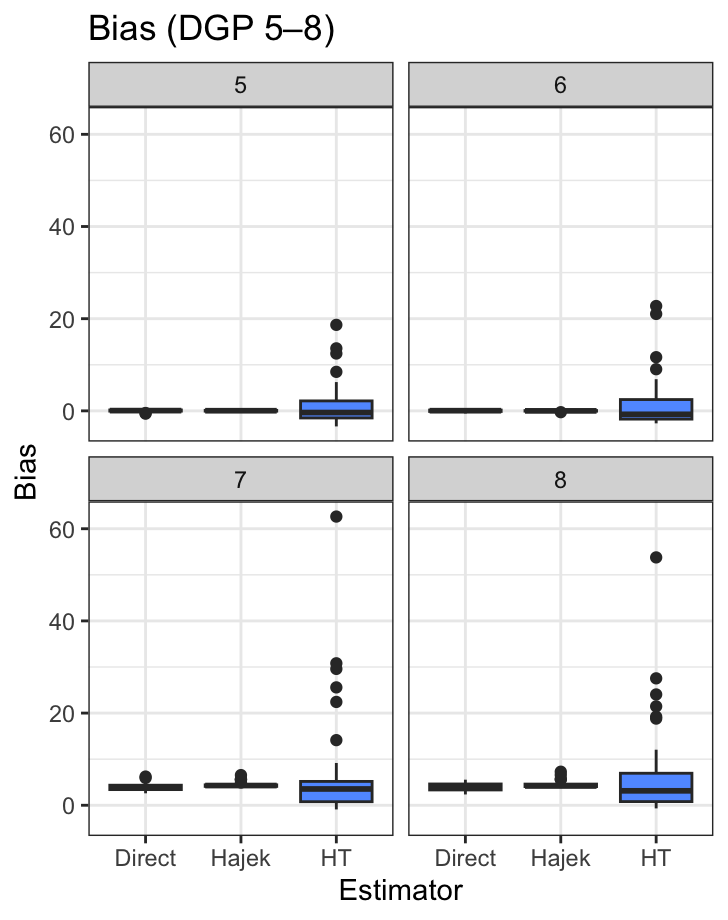}
  \end{minipage}
     \caption{Boxplots of the area-specific estimation bias across the 50 areas for each data generating process (DGP) with total population size $N=400{,}000$.  
For each estimator, the distribution of $\mathrm{Bias}(\hat{\tau}_j)$ across areas is summarized over the Monte Carlo
replications. }
     \label{fig:s1}
\end{figure}
\begin{figure}[t]
  \centering
  \begin{minipage}{0.4\columnwidth}
     \centering
     \includegraphics[width=\columnwidth]{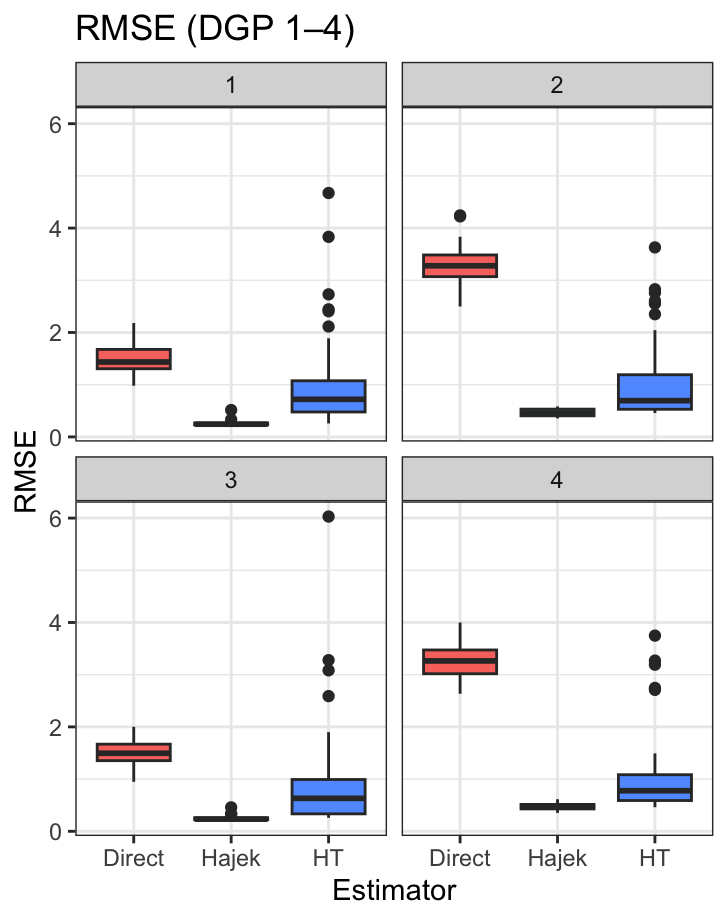}
  \end{minipage}
  \begin{minipage}{0.4\columnwidth}
     \centering
     \includegraphics[width=\columnwidth]{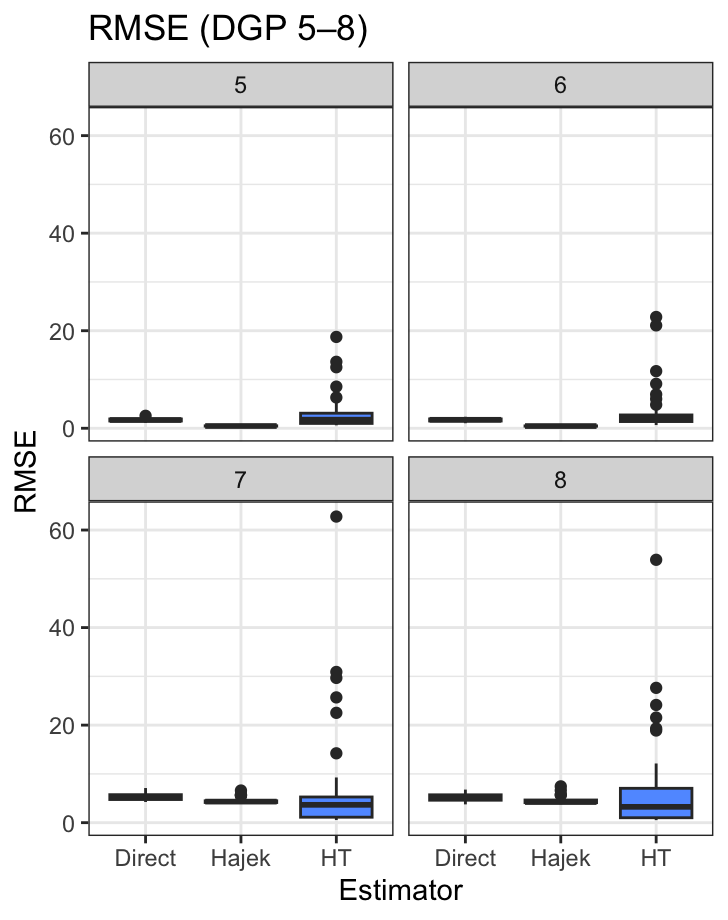}
  \end{minipage}
     \caption{Boxplots of the root mean squared error (RMSE) of the area-specific estimation across the 50 areas for each data generating process (DGP) with total population size $N=400{,}000$. 
For each estimator, the distribution of $\mathrm{RMSE}(\hat{\tau}_j)$ across areas is summarized over the Monte Carlo
replications. }
     \label{fig:s2}
\end{figure}
\begin{figure}[t]
  \centering
  \begin{minipage}{0.4\columnwidth}
     \centering
     \includegraphics[width=\columnwidth]{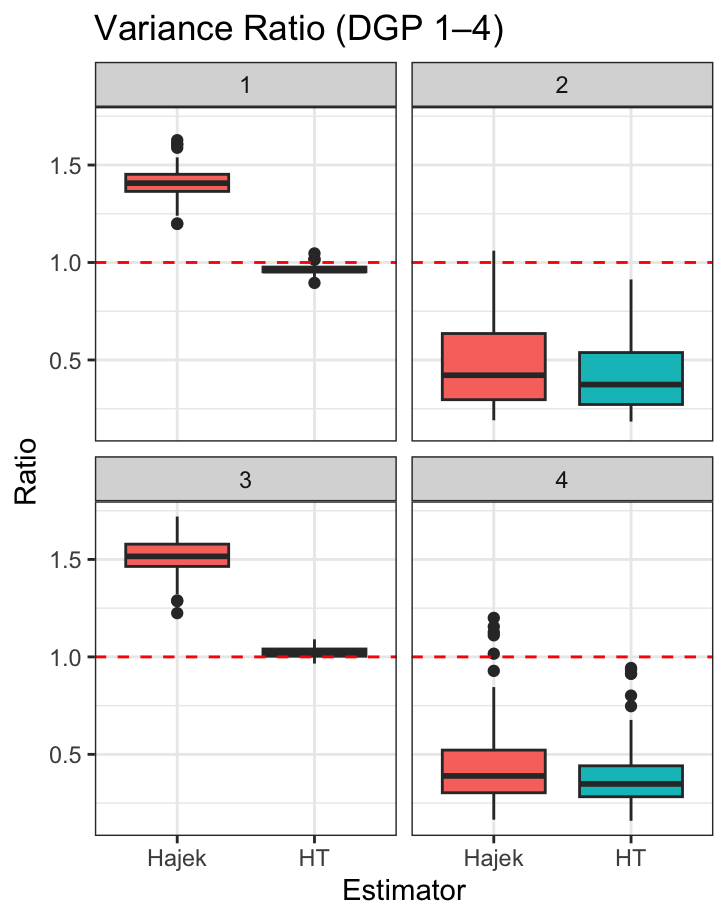}
  \end{minipage}
  \begin{minipage}{0.4\columnwidth}
     \centering
     \includegraphics[width=\columnwidth]{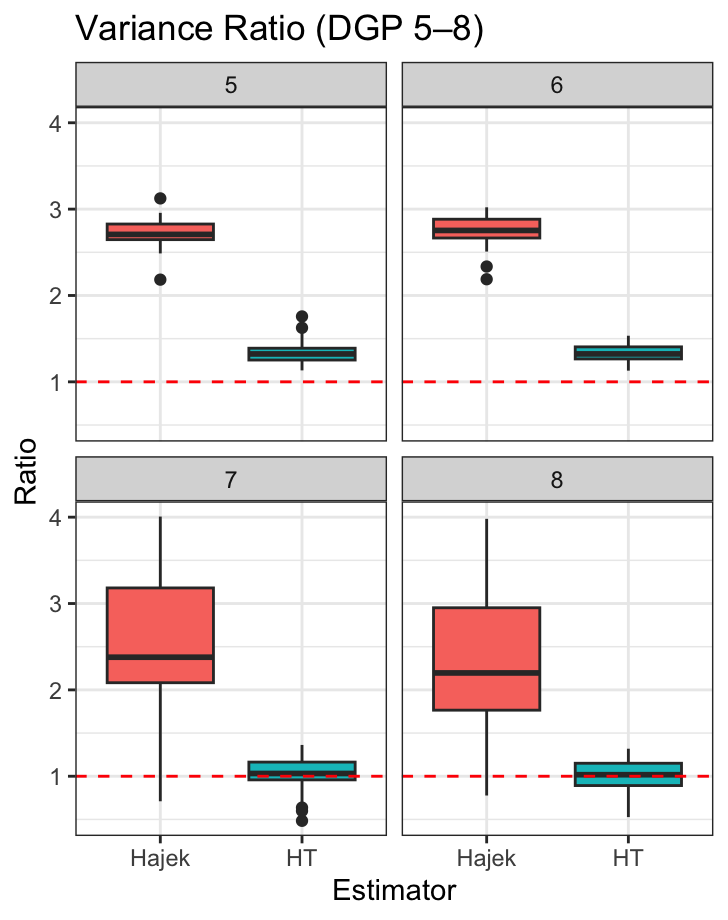}
  \end{minipage}
     \caption{Boxplots of the ratio of variance estimation of the area-specific estimation to the true variance across the 50 areas for each data generating process (DGP) with total population size $N=400{,}000$. 
For each estimator, the distribution of $\mathrm{Ratio}(\hat{V}_j)$ across areas is summarized over the Monte Carlo
replications. }
     \label{fig:s3}
\end{figure}

We focus on the simulation results for the competing estimators under the population size $N=400{,}000$.
Figures \ref{fig:s1}--\ref{fig:s3} respectively show boxplots of the bias, RMSE, and variance estimation accuracy across the 50 areas for each data generating process, and Table \ref{tab:40} summarizes the corresponding averages and standard deviations.
The results for the larger population size $N=800{,}000$ (Figures \ref{fig:s4}--\ref{fig:s6} and Table \ref{tab:80}) are reported in the Supplementary Material. 
The patterns are qualitatively similar across the two settings.

We first examine the bias and RMSE of each estimator under DGP1--DGP5, where the doubly robust condition holds.
In this case, the proposed estimators are expected to be consistent and asymptotically efficient, as established in Theorem \ref{thm:dr} and \ref{thm:seb}. 
However, the left panels of Figure \ref{fig:s1} and Table \ref{tab:40} show that the Horvitz--Thompson type estimator shows substantial variability in bias across areas, leading to unstable performance. 
In fact, its bias is often more dispersed than that of the direct estimator, reflecting the sensitivity of inverse probability weighting to extreme weights in finite samples.
This phenomenon is well documented in the causal inference literature \citep{CCD18}.
In contrast, the proposed H\'ajek type estimator mitigates this instability through normalization, resulting in more stable bias across areas.

With respect to RMSE, the proposed H\'ajek type estimator consistently achieves the smallest RMSE, substantially outperforming both the Horvitz--Thompson type estimator and the direct estimator.
The PRIAL of the H\'ajek type estimator relative to the direct estimator is substantial across all DGP1–DGP5 scenarios, typically ranging from 72.9\% to 86.7\%.
In contrast, the Horvitz--Thompson estimator exhibits highly variable efficiency, with PRIAL ranging from -24.0\% to 76.1\%.
These results indicate that borrowing information across areas can lead to substantial efficiency gains.
Under DGP6--DGP8, where the doubly robust condition is violated, the performance of all estimators deteriorates in terms of both bias and RMSE. 
Nevertheless, the extent of deterioration differs across estimators.
In DGP6, where only the propensity score model is correctly specified, the H\'ajek type estimator does not exhibit a severe degradation in performance and remains relatively stable compared to the other estimators. 
Overall, even in the absence of the doubly robust condition, the normalization in the H\'ajek type estimator appears to mitigate extreme weighting, leading to comparatively stable performance.
In contrast, the Horvitz--Thompson estimator performs poorly, with performance comparable to that of the direct estimator.

Next, we examine the accuracy of the variance estimator for the proposed method, which is shown in Figure \ref{fig:s3}.
The results indicate that the variance estimator performs reasonably well in some cases, particularly under DGP1 and DGP3, where the estimated variance is broadly aligned with the empirical variance.
However, even under correct specification, the performance is not uniformly satisfactory, and the variance estimates based on the H\'ajek type estimator can deviate from the true variance. 
This reflects the additional variability introduced by normalization in the H\'ajek type estimator, which can affect finite-sample variance estimation despite its asymptotic efficiency.
This suggests that variance estimation remains challenging in finite samples, highlighting the need for improved variance estimation methods for reliable inference.

Overall, the simulation results highlight the fundamental limitation of the direct estimator in small-area settings  and demonstrate the advantage of incorporating cross-area information.
The proposed estimator achieves substantial efficiency gains when the doubly robust condition holds, while maintaining relatively stable performance under misspecification. 
At the same time, the results reveal limitations of the Horvitz--Thompson estimator and challenges in variance estimation, which are important for practical applications.

\begin{table}[t]
\caption{Average value of Bias, RMSE and Variance Ratio for each data generating process (DGP) with total population size $N=400{,}000$.
The values in parentheses represent standard deviations.}
\centering

\medskip
\begin{tabular}{
c
S[table-format=1.3] S[table-format=1.3] S[table-format=1.3]
S[table-format=1.3] S[table-format=1.3] S[table-format=1.3]
S[table-format=1.3] S[table-format=1.3]
}
\toprule
\multirow{2}{*}{DGP} & \multicolumn{3}{c}{Bias} & \multicolumn{3}{c}{RMSE} & \multicolumn{2}{c}{Variance Ratio} \\
\cmidrule(lr){2-4} \cmidrule(lr){5-7} \cmidrule(lr){8-9}
& {Direct} & {HT} & {H\'ajek} & {Direct} & {HT} & {H\'ajek} & {HT} & {H\'ajek} \\
\midrule

\multirow{2}{*}{DGP1}
& -0.030 & 0.306 & -0.022 & 2.187 & 1.149 & {\bfseries 0.379} & 0.986 & 1.420 \\
& \multicolumn{1}{c}{\footnotesize (0.211)} & \multicolumn{1}{c}{\footnotesize (1.606)} & \multicolumn{1}{c}{\footnotesize (0.195)}
& \multicolumn{1}{c}{\footnotesize (0.432)} & \multicolumn{1}{c}{\footnotesize (1.243)} & \multicolumn{1}{c}{\footnotesize (0.073)}
& \multicolumn{1}{c}{\footnotesize (0.032)} & \multicolumn{1}{c}{\footnotesize (0.100)} \\
\addlinespace[3pt]

\multirow{2}{*}{DGP2}
& -0.132 & 0.130 & -0.141 & 5.094 & 1.220 & {\bfseries 0.676} & 0.407 & 0.461 \\
& \multicolumn{1}{c}{\footnotesize (0.497)} & \multicolumn{1}{c}{\footnotesize (1.141)} & \multicolumn{1}{c}{\footnotesize (0.155)}
& \multicolumn{1}{c}{\footnotesize (0.684)} & \multicolumn{1}{c}{\footnotesize (0.726)} & \multicolumn{1}{c}{\footnotesize (0.088)}
& \multicolumn{1}{c}{\footnotesize (0.154)} & \multicolumn{1}{c}{\footnotesize (0.189)} \\
\addlinespace[3pt]

\multirow{2}{*}{DGP3}
& 0.031 & 0.335 & 0.009 & 2.197 & 1.213 & {\bfseries 0.370} & 0.975 & 1.439 \\
& \multicolumn{1}{c}{\footnotesize (0.211)} & \multicolumn{1}{c}{\footnotesize (1.932)} & \multicolumn{1}{c}{\footnotesize (0.177)}
& \multicolumn{1}{c}{\footnotesize (0.435)} & \multicolumn{1}{c}{\footnotesize (1.608)} & \multicolumn{1}{c}{\footnotesize (0.056)}
& \multicolumn{1}{c}{\footnotesize (0.027)} & \multicolumn{1}{c}{\footnotesize (0.109)} \\
\addlinespace[3pt]

\multirow{2}{*}{DGP4}
& -0.007 & 0.248 & -0.050 & 5.012 & 1.251 & {\bfseries 0.680} & 0.422 & 0.481 \\
& \multicolumn{1}{c}{\footnotesize (0.539)} & \multicolumn{1}{c}{\footnotesize (1.190)} & \multicolumn{1}{c}{\footnotesize (0.188)}
& \multicolumn{1}{c}{\footnotesize (0.634)} & \multicolumn{1}{c}{\footnotesize (0.780)} & \multicolumn{1}{c}{\footnotesize (0.078)}
& \multicolumn{1}{c}{\footnotesize (0.189)} & \multicolumn{1}{c}{\footnotesize (0.231)} \\
\addlinespace[3pt]

\multirow{2}{*}{DGP5}
& 0.098 & 1.184 & 0.024 & 2.466 & 3.058 & {\bfseries 0.668} & 1.323 & 2.690 \\
& \multicolumn{1}{c}{\footnotesize (0.252)} & \multicolumn{1}{c}{\footnotesize (4.580)} & \multicolumn{1}{c}{\footnotesize (0.182)}
& \multicolumn{1}{c}{\footnotesize (0.425)} & \multicolumn{1}{c}{\footnotesize (3.737)} & \multicolumn{1}{c}{\footnotesize (0.045)}
& \multicolumn{1}{c}{\footnotesize (0.108)} & \multicolumn{1}{c}{\footnotesize (0.169)} \\
\addlinespace[3pt]

\multirow{2}{*}{DGP6}
& 0.071 & 1.077 & -0.044 & 2.520 & 3.115 & {\bfseries 0.678} & 1.335 & 2.698 \\
& \multicolumn{1}{c}{\footnotesize (0.287)} & \multicolumn{1}{c}{\footnotesize (4.393)} & \multicolumn{1}{c}{\footnotesize (0.205)}
& \multicolumn{1}{c}{\footnotesize (0.524)} & \multicolumn{1}{c}{\footnotesize (3.414)} & \multicolumn{1}{c}{\footnotesize (0.083)}
& \multicolumn{1}{c}{\footnotesize (0.097)} & \multicolumn{1}{c}{\footnotesize (0.128)} \\
\addlinespace[3pt]

\multirow{2}{*}{DGP7}
& 4.237 & 5.881 & 4.473 & 6.871 & 6.195 & {\bfseries 4.578} & 0.973 & 2.218 \\
& \multicolumn{1}{c}{\footnotesize (0.888)} & \multicolumn{1}{c}{\footnotesize (6.665)} & \multicolumn{1}{c}{\footnotesize (0.366)}
& \multicolumn{1}{c}{\footnotesize (0.787)} & \multicolumn{1}{c}{\footnotesize (6.579)} & \multicolumn{1}{c}{\footnotesize (0.367)}
& \multicolumn{1}{c}{\footnotesize (0.152)} & \multicolumn{1}{c}{\footnotesize (0.629)} \\
\addlinespace[3pt]

\multirow{2}{*}{DGP8}
& 4.446 & 5.998 & 4.657 & 7.013 & 6.359 & {\bfseries 4.757} & 0.972 & 2.248 \\
& \multicolumn{1}{c}{\footnotesize (0.828)} & \multicolumn{1}{c}{\footnotesize (5.617)} & \multicolumn{1}{c}{\footnotesize (0.401)}
& \multicolumn{1}{c}{\footnotesize (0.804)} & \multicolumn{1}{c}{\footnotesize (5.449)} & \multicolumn{1}{c}{\footnotesize (0.404)}
& \multicolumn{1}{c}{\footnotesize (0.151)} & \multicolumn{1}{c}{\footnotesize (0.657)} \\

\bottomrule
\end{tabular}
\label{tab:40}
\end{table}

\section{Application to 2024 presidential election data}\label{sec:rda}

To illustrate the practical implications of the the proposed method, we analyze survey data from the 2024 wave of the American National Election Studies (ANES) \citep{ANES}. 
The ANES is a nationally representative survey of the U.S. voting-age population that collects detailed information on political attitudes, campaign experiences and demographic characteristics.
Estimating the effect of campaign contact on voters' candidate evaluations at the state level is of substantive interest. However, both campaign strategies and voter responses are likely to vary across states, reflecting differences in political environments and electoral competitiveness. 
As a result, the treatment effect is expected to exhibit substantial heterogeneity across states.

Despite the rich individual-level information in the ANES, the sample size within each state is limited. 
Consequently, direct estimation of state-level treatment effects is often unstable, and for most states infeasible due to the absence of treated or control observations or lack of overlap in propensity scores.
To address these challenges, we apply the proposed small area estimation approach, which borrows strength across states while allowing for heterogeneity in treatment effects.

\subsection{Data description}
The outcome variable is based on respondents' feeling thermometer ratings evaluated on a 0--100 scale toward the two major presidential candidates in the 2024 election.
We construct the outcome as the difference between the assessments of Donald Trump and that of Kamala Harris, $Y_i = FT_i^{\rm Trump} - FT_i^{\rm Harris}$, which takes values between $-100$ and $100$.
Positive values indicate warmer feelings toward Trump relative to Harris, while negative values indicate warmer feelings toward Harris.
The treatment variable $T_i$ indicates whether the respondent was contacted by a political campaign from both the Democratic Party and the Republican Party during the election cycle.

We consider two sets of covariates.
The first set, denoted by $\X_i$, consists of demographic variables that are available both in the survey and from external population data sources such as the U.S. Census. 
These variables include age, sex, race and ethnicity, educational attainment and total household income.
The second set of covariates, $\Z_i$, includes seven survey-based measures of political predispositions and behavior not available in census data: party identification, ideological self-placement, attention to politics, interest in following political campaigns, political knowledge, past political participation and trust in news media.
The ANES 2024 dataset contains $5{,}521$ respondents.
After excluding individuals with missing values in the outcome, treatment, or covariates, the analytical sample consists of $n = 3{,}100$ observations.
Across the 51 states, the average number of observations per state and its standard deviation are respectively 54.0 and 47.4 for the control group and 6.75 and 7.59 for the treated group, indicating substantial imbalance in sample sizes across treatment status.
More details of these variables and their summary statistics are provided in the Supplementary Material.

\subsection{Estimation procedure}

We estimate the state-level treatment effects using the H\'ajek type estimator defined in (\ref{eqn:HT}). 
The implementation closely follows the simulation design in Section \ref{sec:sim}, with the key difference that the sampling probability is unknown in the empirical application and must be approximated using survey weights.
Specifically, we estimate the outcome regression functions, the propensity score and the area assignment probabilities using standard parametric models (linear, logistic, and multinomial logistic), all based on the survey sample.
To approximate the sampling probability, we decompose $\mathbb{P}(S=1 | \X,A=j)$ as
$$
\mathbb{P}(S=1 | \X,A=j) = \frac{\mathbb{P}(S=1 | \X)\, \mathbb{P}(A=j | \X,S=1)}{\mathbb{P}(A=j | \X)}.
$$
We proxy $\mathbb{P}(S=1 | \X)$ using the inverse of the ANES raked survey weights, following standard practice in survey sampling \citep{S03}.
The conditional probability $\mathbb{P}(A=j | \X,S=1)$ is estimated from the survey via a multinomial logit model. 
The marginal distribution $\mathbb{P}(A=j | \X)$ is approximated using population data from the American Community Survey (ACS) \citep{ACS} by partitioning $\X$ into discrete cells and computing the corresponding population shares.

The estimated components are substituted into the sample analogue of the doubly robust score function. 
We further trim observations with extreme values of the estimated doubly robust score excluding those exceeding $800$ in absolute value.
Following this procedure, the resulting sample sizes are highly comparable across states, ranging from $3{,}002$ to $3{,}100$.
We then compute the H\'ajek type estimator of each state in \eqref{eqn:HT} and its standard error based on variance  estimator in \eqref{eqn:vht}.

\subsection{Assessment of the area ignorability assumption}

\begin{figure}[h]
\centering
\begin{minipage}[b]{0.49\columnwidth}
    \centering
    \includegraphics[width=0.9\columnwidth]{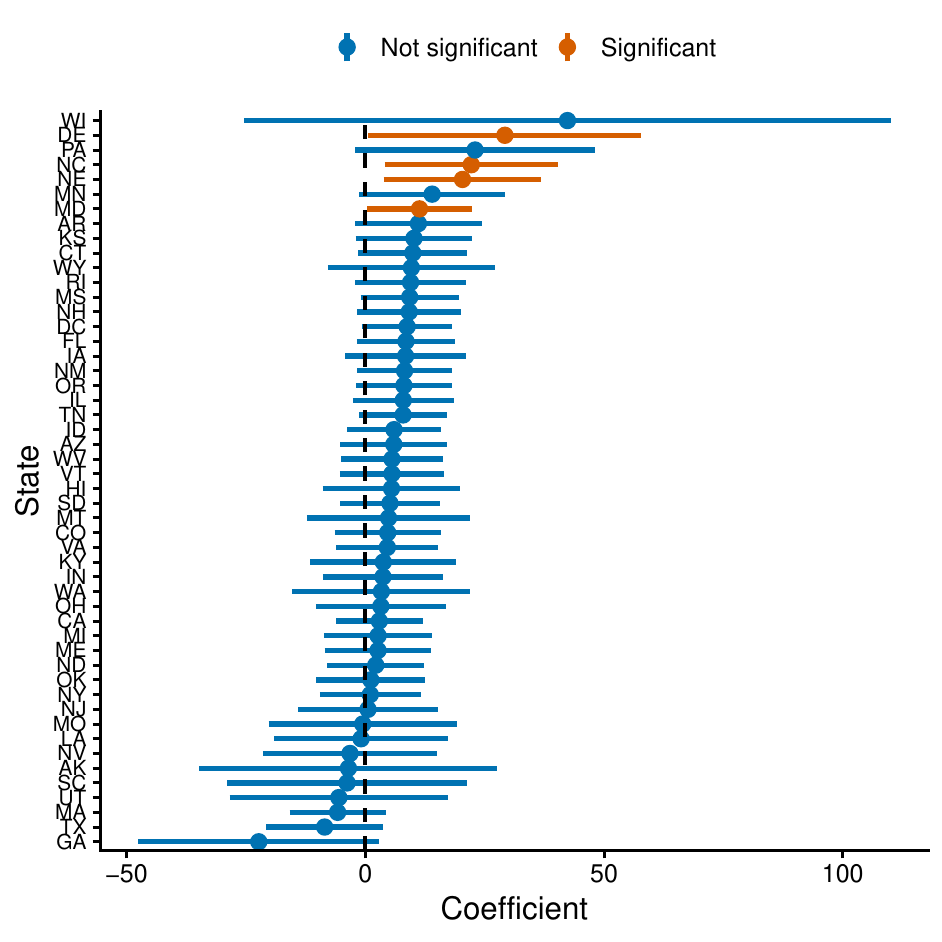}
\end{minipage}
\begin{minipage}[b]{0.49\columnwidth}
    \centering
    \includegraphics[width=0.9\columnwidth]{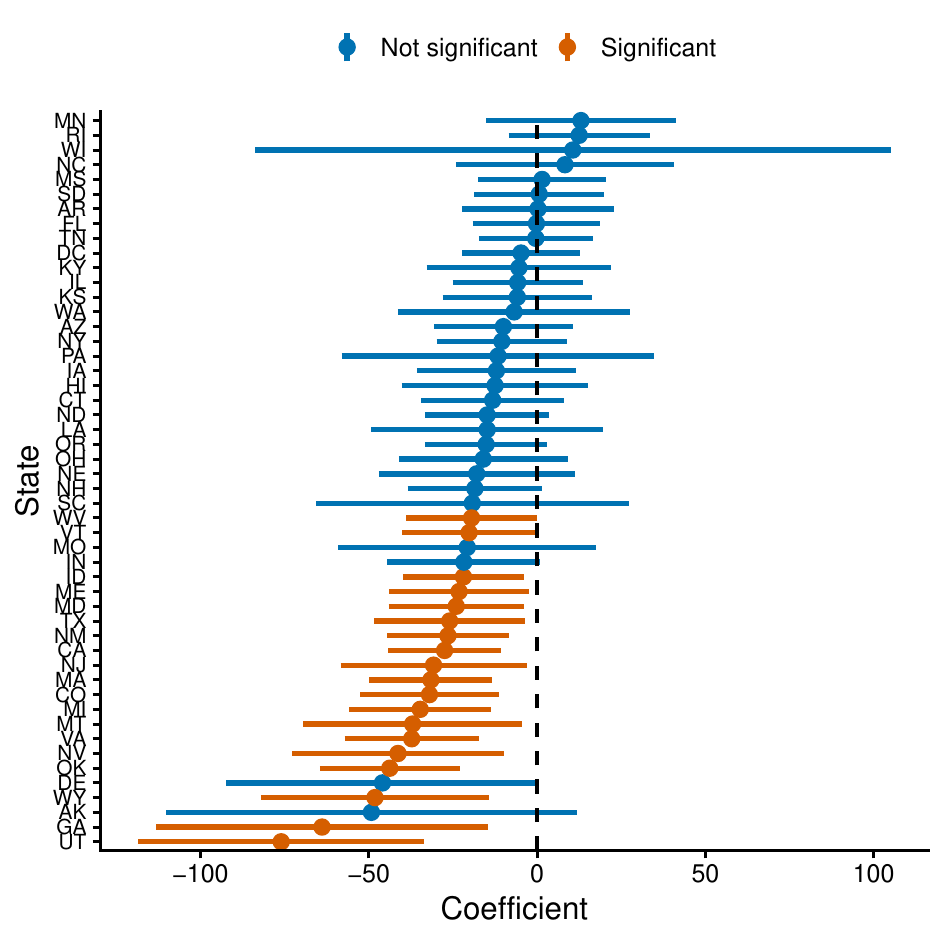}
\end{minipage}
    \caption{Estimated coefficients for state indicators from a regression of the outcome $Y$ on state indicators $A$, demographic covariates $\X$, and political covariates $\Z$ (left figure) and from a regression of the outcome $Y$ on state indicators $A$ and demographic covariates $\X$ (right figure).
Points represent coefficient estimates and horizontal lines indicate 95\% confidence intervals.}
    \label{fig:aaia}
\end{figure}

Our identification strategy relies on the assumption that the potential outcomes are conditionally independent of the area indicator given $(\X,\Z,S=1)$, namely $Y(t) \perp A | \X,\Z,S=1$ for $t=0,1$.
To examine the plausibility of this assumption, we conduct a simple diagnostic analysis based on regression models, similar to that conducted in \citet{KY21}.
The left figure in Figure~\ref{fig:aaia} presents the 95\% confidence intervals of the coefficients for the area indicators obtained from a regression of the outcome $Y$ on the area indicators $A$, the demographic covariates $\X$, and the political covariates $\Z$, using only sampled units ($S=1$). 
The coefficients correspond to the indicators for $A=2,\ldots,51$, with the first area serving as the reference category.
For comparison, the right figure in Figure \ref{fig:aaia} reports the corresponding confidence intervals from a regression that includes only the area indicators and the demographic covariates $\X$, but excludes the political covariates $\Z$.

The comparison between two figures in Figure \ref{fig:aaia} suggests that including the political covariates substantially reduces systematic differences across areas. 
When conditioning only on $\X$, many of the area coefficients appear significantly different from zero, indicating residual cross-state differences in the outcome. 
However, once the additional political covariates $\Z$ are included, most of the area coefficients become statistically indistinguishable from zero.
This pattern suggests that differences in political predispositions across states largely account for the observed cross-state variation in the outcome. 
Therefore, the assumption that the potential outcomes are conditionally independent of the area indicator given $(\X,\Z,S=1)$ appears reasonably plausible in this application.

\subsection{Estimated area-specific treatment effects}

\begin{figure}[h]
\centering
\begin{minipage}[b]{0.49\columnwidth}
    \centering
    \includegraphics[width=1.0\columnwidth]{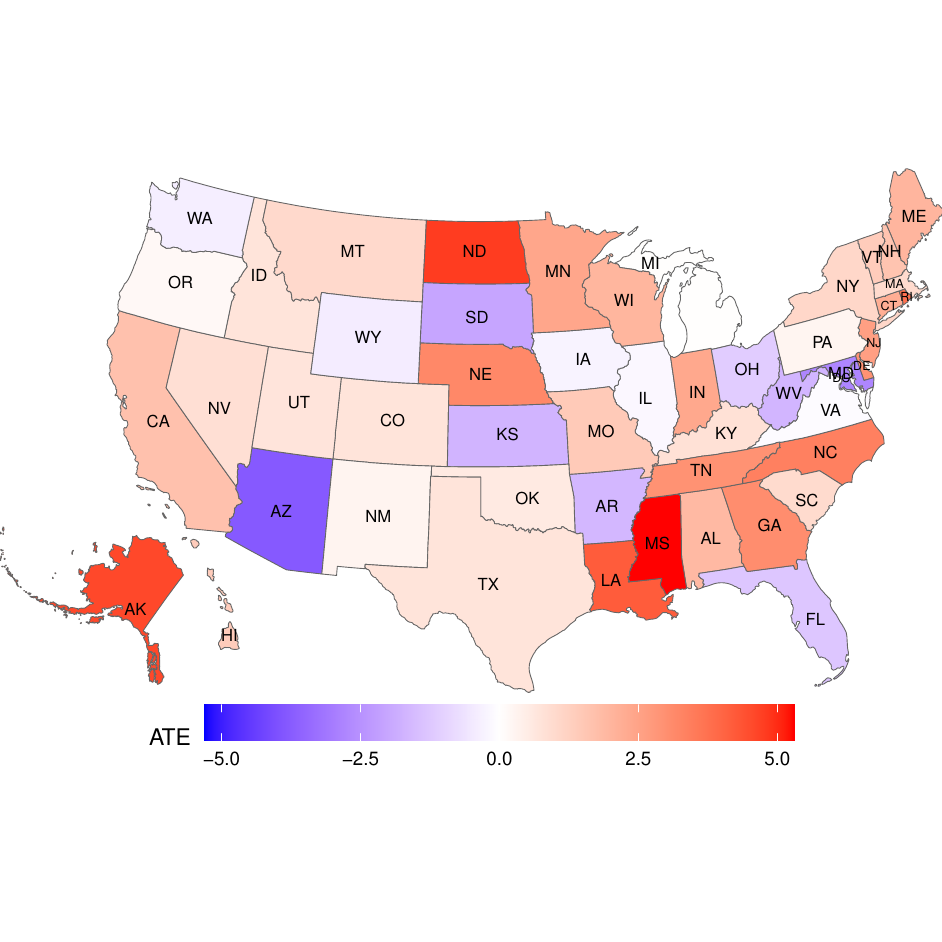}
\end{minipage}
\begin{minipage}[b]{0.49\columnwidth}
    \centering
    \includegraphics[width=1.0\columnwidth]{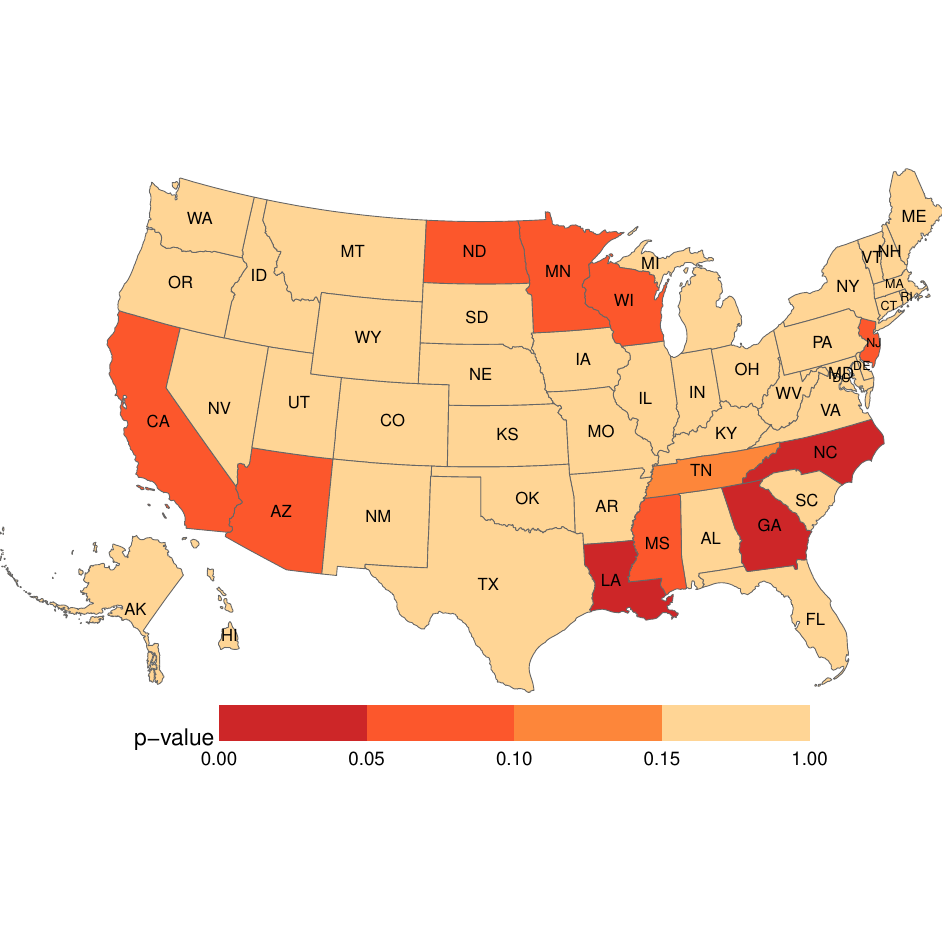}
\end{minipage}
    \caption{State-level average treatment effects (ATEs) of campaign contact and their statistical significance. 
The left panel reports the estimated ATEs for each state, while the right panel displays the corresponding $p$-values. }
    \label{fig:map}
\end{figure}

Figure \ref{fig:map}  presents the estimated state-specific average treatment effects of campaign contact on the feeling thermometer difference between Trump and Harris along with their corresponding $p$-values.
Overall, the results suggest that the causal effect of campaign contact is limited across most states, which is consistent with a large body of empirical findings in political science showing modest persuasion effects \citep{KB18}.
Despite this overall null pattern, the results reveal substantial heterogeneity across states. 
Notably, among the seven battleground states, four states—Arizona, Georgia, North Carolina, and Wisconsin—exhibit statistically significant treatment effects.

Furthermore, as in Table \ref{tab:dp}, to benchmark the proposed estimator, we construct direct estimates of the state-specific ATE for five states where the direct estimate can be constructed.
The direct estimates exhibit extremely large standard errors, reflecting instability of within-area estimation.
In contrast, the proposed estimator yields substantially smaller standard errors, indicating efficiency gains from borrowing information across areas.
Moreover, while the direct estimates fluctuate widely and take implausibly large values, the proposed estimates are more stable and lie within a reasonable range, suggesting improved reliability.
Overall, the results highlight the limitations of direct estimation in small-area settings and demonstrate the practical advantages of the proposed approach.

\begin{table}[h]
\centering
\caption{The direct estimates and the estimates obtained by the proposed method and their standard deviations, along with the sample sizes for treated ($T=1$) and control ($T=0$) groups for the five selected states.}
\begin{tabular}{lccccccc}
\toprule
\multirow{2}{*}{State} & \multicolumn{2}{c}{Sample size} & \multicolumn{2}{c}{Direct estimate} & \multicolumn{2}{c}{Proposed method} \\
& $T=1$ & $T=0$ & estimate & s.d. & estimate & s.d. \\
\midrule
California & 25 & 223 & -4.21 & 130.03 & 1.71 & 1.04 \\
Massachusetts & 28 & 91 & 0.90 & 40.49 & 1.00 & 2.19 \\
New York & 18 & 66 & -5.32 & 174.20 & 1.10 & 1.27 \\
Oregon & 25 & 95 & 6.56 & 23.74 & -0.19 & 2.29 \\
West Virginia & 29 & 64 & -8.23 & 19.20 & -1.68 & 3.08  \\
\bottomrule
\end{tabular}
\label{tab:dp}
\end{table}

\section{Concluding Remarks}\label{sec:con}

In this paper, we develop a new framework for causal small area estimation in settings where both treatment and outcome are observed only for sampled units, while auxiliary covariates are available at the population level.
By integrating design-based inference with causal inference methods, we establish identification of area-specific treatment effects and propose a doubly robust estimator that remains consistent under partial misspecification of nuisance components.
The proposed approach enables the use of rich survey-only covariates, thereby extending the scope of small area estimation beyond traditional settings.
We derive the semiparametric efficiency bound for the target parameter and show that the proposed estimator attains this bound under regularity conditions.
The use of cross-fitting and orthogonal scores allows for flexible estimation of nuisance functions while preserving valid asymptotic inference.
Together, these results provide a principled foundation for causal inference in small area settings under realistic data constraints.

An important direction for future research is the development of accurate variance estimators for the proposed procedure.
Constructing inference methods that remain reliable in finite samples and properly account for complex sampling designs and nuisance estimation is an important open problem.
Moreover, when the number of areas is large, the area assignment probability requires a large number of parameters, which may lead to overfitting of the model. 
To mitigate the issue, it would be an interesting future work to develop parsimonious modeling.

\section*{Acknowledgement}
This paper is partially supported by the Japan Society of the Promotion of Science (JPSP KAKENHI) grant numbers 25H00546 and 25K16613.

\bibliographystyle{chicago}
\bibliography{ref}

\newpage


\vspace{0.5cm}
\begin{center}
{\LARGE {\bf Supplementary Material for ``Causal Small Area Estimation with Survey-only Covariates"}}
\end{center}

\setcounter{page}{1}
\setcounter{equation}{0}
\renewcommand{\theequation}{S\arabic{equation}}
\setcounter{section}{0}
\renewcommand{\thesection}{S\arabic{section}}
\setcounter{lem}{0}
\renewcommand{\thelem}{S\arabic{lem}}
 \setcounter{table}{0}
\renewcommand{\thetable}{S\arabic{table}}
\setcounter{figure}{0}
\renewcommand{\thefigure}{S\arabic{figure}}

\section{Proof of Proposition \ref{prp:id}}

The estimand of the target area $j$ can be written as 
\begin{align*}
&\tau_j=
\mathbb{\mathbb{E}}[Y(1)-Y(0)|A=j] \\
& =
\mathbb{\mathbb{E}}[\mathbb{\mathbb{E}}[Y(1)-Y(0)|\X,S=1,A=j] |A=j] \\
& =
\mathbb{\mathbb{E}}[\mathbb{\mathbb{E}}[\mathbb{\mathbb{E}}[Y(1)-Y(0) |\X,\Z,S=1] |\X,S=1,A=j] |A=j],
\end{align*}
where the second and third equalities hold from Assumption 1 (Sampling ignorability) and Assumption 4 (Area ignorability), respectively.

From the standard theory of causal inference, under Assumption 2 (Unconfoundedness) and Assumption 3 (Overlap), $\mathbb{\mathbb{E}}[Y(1)-Y(0) |\X,\Z,S=1]$ can be identified as
\begin{align*}
\mathbb{\mathbb{E}}[Y(1)-Y(0) |\X,\Z,S=1]
=
\mathbb{\mathbb{E}}\bigg[\frac{TY}{e(\X,\Z)}-\frac{(1-T)Y}{1-e(\X,\Z)} \bigg|\X,\Z,S=1\bigg] 
\end{align*}
Then, by following the proof of Proposition 1 in Kuriwaki and Yamauchi (2021), the proposition holds.

\section{Proof of Theorem \ref{thm:dr}}

Remind that for $\widetilde{\bta}=(\widetilde{m}_1,\widetilde{m}_0,\widetilde{e},\widetilde{\pi}_A)$, 
\begin{align}\label{eqn:tr}
&{\mathbb{E}}[\phi_{j}(\W, \widetilde{\bta})] \nonumber \\ 
\begin{split}
&=
\mathbb{\mathbb{E}}\bigg[\frac{\one\{S=1\}}{\mathbb{P}(A=j)} \frac{\widetilde{\pi}_A(j | \X,\Z)}{\mathbb{P}(A=j|\X,S=1)}\frac{\mathbb{P}(A=j|\X)}{\mathbb{P}(S=1|\X)} \\
&\times
\bigg\{\frac{T\{Y-\widetilde{m}_1(\X,\Z)\}}{\widetilde{e}(\X,\Z)}
-\frac{(1-T)\{Y-\widetilde{m}_0(\X,\Z)\}}{1-\widetilde{e}(\X,\Z)}
+\widetilde{m}_1(\X,\Z) -\widetilde{m}_0(\X,\Z)\bigg\}\bigg] \\
& \quad +\mathbb{\mathbb{E}}\bigg[\frac{\one\{S=1\}}{\mathbb{P}(A=j)} \frac{\one\{A=j\}-\widetilde{\pi}_A(j | \X,\Z)}{\mathbb{P}(A=j|\X,S=1)}\frac{\mathbb{P}(A=j|\X)}{\mathbb{P}(S=1|\X)} 
\{\widetilde{m}_1(\X,\Z) -\widetilde{m}_0(\X,\Z)\} \bigg].
\end{split}
\end{align}

First, we consider the case, in which the area assignment probability is correctly specified, namely, $\widetilde{\pi}_A(j | \X,\Z)=\pi_A(j|\X,\Z)=\mathbb{P}(A=j | \X, \Z, S=1)$ holds.
Then, we have
\begin{align*}
&\mathbb{\mathbb{E}}[\one\{S=1\}\{\one\{A=j\}-\pi_A(j\X,\Z)\}|\X,\Z] \\
&=
\mathbb{P}(S=1,A=j|\X,\Z)-\mathbb{P}(S=1|\X,\Z)\mathbb{P}(A=j | \X, \Z, S=1)
=0,
\end{align*}
which implies that the second term in \eqref{eqn:tr} vanishes.
Additionally, if the outcome regression model is correctly specified, but the propensity score is misspecified, namely $\widetilde{m}_t=m_t$ for $t=0,1$ and $\widetilde{e}=e^*\neq e$ hold, we have
\begin{align*}
&{\mathbb{E}}\bigg[\frac{T\{Y-m_1(\X,\Z)\}}{e^*(\X,\Z)}
+m_1(\X,\Z) \bigg|\X,\Z,S=1\bigg] \\
&=
{\mathbb{E}}\bigg[\frac{T\{{\mathbb{E}}[Y|\X,\Z,S=1,T=1]-m_1(\X,\Z)\}}{e^*(\X,\Z)}\bigg|\X,\Z,S=1\bigg]
+{\mathbb{E}}[Y(1)|\X,\Z,S=1] \\
&=
{\mathbb{E}}[Y(1)|\X,\Z,S=1].
\end{align*}
Similarly, we have
\begin{align*}
\mathbb{\mathbb{E}}\bigg[\frac{(1-T)\{Y-m_0^S(\X,\Z)\}}{1-e^*(\X,\Z)}
+m_0^S(\X,\Z)\bigg]
=
\mathbb{\mathbb{E}}[Y(0)|\X,\Z,S=1].
\end{align*}
Combining the above two results, we have the following result.
\begin{align*}
&\mathbb{\mathbb{E}}\bigg[\frac{T\{Y-m_1^S(\X,\Z)\}}{e^*(\X,\Z)}
-\frac{(1-T)\{Y-m_0^S(\X,\Z)\}}{1-e^*(\X,\Z)}
+m_1^S(\X,\Z) -m_0^S(\X,\Z)\bigg|\X,\Z,S=1\bigg] \\
&=
\mathbb{E}[Y(1)-Y(0)|\X,\Z,S=1].
\end{align*}
Then, the right-hand side of \eqref{eqn:tr} can be written as
\begin{align*}
&
\mathbb{E}\bigg[\frac{\one\{S=1\}}{\mathbb{P}(A=j)} \frac{\pi_A(j | \X,\Z)}{\mathbb{P}(A=j|\X,S=1)}\frac{\mathbb{P}(A=j|\X)}{\mathbb{P}(S=1|\X)} 
\mathbb{E}[Y(1)-Y(0)|\X,\Z,S=1]\bigg] \\
&=
\mathbb{E}[Y(1)-Y(0)|A=j] = \tau_j,
\end{align*}
where the first equality holds from the proof of Proposition 1 in Kuriwaki and Yamauchi (2021).

On the other hand, if the propensity score is correctly specified but the outcome regression model is misspecified, that is, $\widetilde{e}=e$ and $\widetilde{m}_t= m_t^* \neq m_t$ for $t=0,1$, 
\begin{align*}
&\mathbb{E}\bigg[\frac{T\{Y-m_1^*(\X,\Z)\}}{e(\X,\Z)}
+m_1^*(\X,\Z) \bigg|\X,\Z,S=1\bigg] \\
&=
\mathbb{E}\bigg[\frac{TY}{e(\X,\Z)}
+\bigg(1-\frac{T}{e(\X,\Z)}\bigg)m_1^*(\X,\Z) \bigg|\X,\Z,S=1\bigg] \\
&=
\mathbb{E}\bigg[\frac{TY}{e(\X,\Z)}\bigg|\X,\Z,S=1\bigg] 
+
\bigg(1-\frac{\mathbb{E}[T|\X,\Z,S=1]}{e(\X,\Z)}\bigg)m_1^*(\X,\Z) \\
&=
\mathbb{E}[Y(1)|\X,\Z,S=1].
\end{align*}
Similarly, we have
\begin{align*}
\mathbb{E}\bigg[\frac{(1-T)\{Y-m_0^*(\X,\Z)\}}{e(\X,\Z)}
+m_0^*(\X,\Z) \bigg|\X,\Z,S=1\bigg] 
=
\mathbb{E}[Y(0)|\X,\Z,S=1].
\end{align*}
Then, from the same argument as above, the right-hand side of \eqref{eqn:tr} coincides with $\tau_j$.

Next, we assume that the area assignment probability is misspecified, but the outcome regression model is correctly specified, namely $\widetilde{\pi}_A=\pi_A^*\neq\pi_A$ and $\widetilde{m}_t=m_t$ for $t=0,1$.
In this case, the right-hand side of \eqref{eqn:tr} can be rewritten as
\begin{align*}
&
\mathbb{E}\bigg[\frac{\one\{S=1\}}{\mathbb{P}(A=j)} \frac{\pi_j^*(\X,\Z)}{\mathbb{P}(A=j|\X,S=1)}\frac{\mathbb{P}(A=j|\X)}{\mathbb{P}(S=1|\X)} \\
& \qquad \times
\bigg\{\frac{T\{Y-m_1(\X,\Z)\}}{\widetilde{e}(\X,\Z)}
-\frac{(1-T)\{Y-m_0(\X,\Z)\}}{1-\widetilde{e}(\X,\Z)}\bigg\}\bigg] \\
& +\mathbb{E}\bigg[\frac{\one\{S=1\}}{\mathbb{P}(A=j)} \frac{\one\{A=j\}}{\mathbb{P}(A=j|\X,S=1)}\frac{\mathbb{P}(A=j|\X)}{\mathbb{P}(S=1|\X)} 
\{m_1(\X,\Z) -m_0(\X,\Z)\} \bigg].
\end{align*}
The first term vanishes since the outcome regression models are correctly specified.
By following the first equation in the proof of Proposition 1 in Kuriwaki and Yamauchi (2021), the second term can be written as
\begin{align*}
&\mathbb{E}\bigg[\frac{\one\{S=1\}}{\mathbb{P}(A=j)} \frac{\one\{A=j\}}{\mathbb{P}(A=j|\X,S=1)}\frac{\mathbb{P}(A=j|\X)}{\mathbb{P}(S=1|\X)} 
\{m_1(\X,\Z) -m_0(\X,\Z)\} \bigg] \\
&=
\mathbb{E}\bigg[\frac{\mathbb{P}(A=j|\X)}{\mathbb{P}(A=j)} \frac{\one\{S=1,A=j\}}{\mathbb{P}(S=1,A=j|\X)}
\mathbb{E}[Y(1)-Y(0)|\X,\Z,S=1]\bigg] \\
&=
\mathbb{E}[Y(1)-Y(0)|A=j]=\tau_j.
\end{align*}

Then, under the condition in Theorem \ref{thm:dr}, we have $\tau_j={\mathbb{E}}[\phi_{j}(\W, \widetilde{\bta})]$.

\section{Proof of Theorem \ref{thm:seb}}

From Assumptions (A1) and (A4), we can denote the joint distribution of $\W=(Y,\X,\Z,T,A,S)$ as
$$
p(\W)
=
p_Y(y|\x,\z,t,s)
p_T(t|\x,\z,s)
p_A(a|\x,\z,s)
q(s,\x,\z).
$$
Consider a regular parametric submodel $p_\epsilon$ with score
$
s(\W)=s_Y+s_T+s_A+s_q
$.
The nuisance tangent space of this model is 
$
\mathcal{T}
=
\mathcal{T}_Y \oplus \mathcal{T}_T \oplus \mathcal{T}_A \oplus \mathcal{T}_q,
$
where
\begin{align*}
&\mathcal{T}_Y = \{s_Y: \mathbb{E}[s_Y | \X,\Z,T,S]=0\}, \quad
\mathcal{T}_T = \{s_T: \mathbb{E}[s_T | \X,\Z,S]=0\}, \\
&\mathcal{T}_A = \{s_A: \mathbb{E}[s_A | \X,\Z,S]=0\}, \quad
\mathcal{T}_q = \{s_q: \mathbb{E}[s_q]=0\}.
\end{align*}
From Proposition \ref{prp:id}, the parameter $\tau_j$ is given by
$
\tau_j
=
\mathbb{E}\left[
G(\W)H(\W)
\right],
$
where we define
$$
G(\W)
=
\frac{\one\{S=1\}}{\pi_S(\X,A=j)} \frac{p_A}{p(j)},
\quad
H(\W)
=
\frac{TY}{p_T} - \frac{(1-T)Y}{1-p_T},
$$
and $\pi_S(\X,A=j)$ and $p(j)$ are assumed to be known.
We will construct $\varphi(\W)$, which satisfies pathwise differentiability and orthogonality.

First, we consider the contribution of $s_Y$.
Since for $m_t=m_t(\X,\Z)=\mathbb{E}[Y(t)|\X,\Z,S=1]$  
$$
\mathbb{E}\left[G
\left\{
\frac{Tm_1}{p_T}
-
\frac{(1-T)m_0}{1-p_T}
\right\} s_Y \right] 
=
\mathbb{E}\left[G
\left\{
\frac{Tm_1}{p_T}
-
\frac{(1-T)m_0}{1-p_T}
\right\} \mathbb{E}[s_Y|\X,\Z,T,S] \right] 
=
0,
$$
it holds that
$$
\dot{\tau}_{j(Y)} = \mathbb{E}\left[G
\left\{
\frac{TY}{p_T}
-
\frac{(1-T)Y}{1-p_T}
\right\} s_Y \right] 
=
\mathbb{E}\left[G
\left\{
\frac{T(Y-m_1)}{p_T}
-
\frac{(1-T)(Y-m_0)}{1-p_T}
\right\} s_Y \right]. 
$$
Define
$$
\varphi_{Y} = \frac{\one\{S=1\}}{\pi_S}
\frac{p_A}{p(j)}
\left\{
\frac{T(Y-m_1)}{p_T}
-
\frac{(1-T)(Y-m_0)}{1-p_T}
\right\} .
$$
Since $s_T$ does not depend on $Y$, from the property of score
$$
\mathbb{E}[\varphi_Y s_T]
=
\mathbb{E}\Big[
\mathbb{E}[\varphi_Y s_T | \X,\Z,T,S]
\Big]
=
\mathbb{E}\Big[s_T
\mathbb{E}[\varphi_Y  | \X,\Z,T,S]
\Big]
=0.
$$
Since $\mathbb{E}[\varphi_Y | \X,\Z,S]=0$, we have
$$
\mathbb{E}[\varphi_Y s_A]
=
\mathbb{E}\Big[
\mathbb{E}[\varphi_Y s_A | \X,\Z,S]
\Big]
=
\mathbb{E}\Big[s_A
\mathbb{E}[\varphi_Y | \X,\Z,S]
\Big]
=0,
$$
from Assumption 4 (Area ignorability), and
$
\mathbb{E}[\varphi_Y s_q]=0
$.

Next, we consider the contribution of $s_T$, which appears in $\mathbb{E}[G{\dot H}_T]$ and $\mathbb{E}[GHs_T]$.
Since we have
$$
{\dot H}_{T} = 
-\frac{TY}{p_T^2(\X,\Z)}\dot{p}_T-\frac{(1-T)Y}{(1-p_T)^2}\dot{p}_T
$$
and $\dot{p}_T=\mathbb{E}[(T-p_T)s_T|\X,\Z,S]$ from a property of score function, it holds that
\begin{align*}
&\mathbb{E}[G {\dot H}_T] 
=
\mathbb{E}\left[G \left\{-\frac{TY}{p_T^2(\X,\Z)}\dot{p}_T-\frac{(1-T)Y}{(1-p_T)^2}\dot{p}_T\right\} \right] \\
&=
-\mathbb{E}\left[G \frac{TY}{p_T^2}\mathbb{E}[(T-p_T)s_T|\X,\Z,S]\right] 
-\mathbb{E}\left[G \frac{(1-T)Y}{(1-p_T)^2}\mathbb{E}[(T-p_T)s_T|\X,\Z,S] \right] \\
&=
-\mathbb{E}\left[\mathbb{E}\left[\frac{\mathbb{P}(S=1|\X,\Z)}{\pi_S(\X,A=j)} \frac{p_A}{p(j)} \frac{TY}{p_T^2}\mathbb{E}[(T-p_T)s_T|\X,\Z,S]\right]\bigg|\X,\Z,S=1\right] \\
&\quad -\mathbb{E}\left[\mathbb{E}\left[\frac{\mathbb{P}(S=1|\X,\Z)}{\pi_S(\X,A=j)} \frac{p_A}{p(j)} \frac{(1-T)Y}{(1-p_T)^2}\mathbb{E}[(T-p_T)s_T|\X,\Z,S] \right]\bigg|\X,\Z,S=1\right],
\end{align*}
where the last equality holds from iterated expectation.
By using Lemma A.1 in Kuriwaki and Yamauchi (2021), the last expression of the above equation is
\begin{align*}
&
-\mathbb{E}\left[\mathbb{E}\left[\frac{\mathbb{P}(S=1|\X,\Z)}{\pi_S(\X,A=j)} \frac{p_A}{p(j)} \mathbb{E}\left[\frac{TY}{p_T^2}|\X,\Z,S=1\right](T-p_T)s_T\right]\bigg|\X,\Z,S=1\right] \\
&\quad -\mathbb{E}\left[\mathbb{E}\left[\frac{\mathbb{P}(S=1|\X,\Z)}{\pi_S(\X,A=j)} \frac{p_A}{p(j)} \mathbb{E}\left[\frac{(1-T)Y}{(1-p_T)^2}\bigg|\X,\Z,S=1\right](T-p_T)s_T \right]\bigg|\X,\Z,S=1\right] \\
&=
-\mathbb{E}\left[\mathbb{E}\left[\frac{\mathbb{P}(S=1|\X,\Z)}{\pi_S(\X,A=j)} \frac{p_A}{p(j)} \frac{m_1}{p_T}(T-p_T)s_T\right]\bigg|\X,\Z,S=1\right] \\
&\quad -\mathbb{E}\left[\mathbb{E}\left[\frac{\mathbb{P}(S=1|\X,\Z)}{\pi_S(\X,A=j)} \frac{p_A}{p(j)} \frac{m_0}{(1-p_T)}(T-p_T)s_T \right]\bigg|\X,\Z,S=1\right] \\
&=
-\mathbb{E}\left[\frac{\one\{S=1\}}{\pi_S(\X,A=j)} \frac{p_A}{p(j)} \frac{m_1}{p_T}(T-p_T)s_T\right] 
-\mathbb{E}\left[\frac{\one\{S=1\}}{\pi_S(\X,A=j)} \frac{p_A}{p(j)} \frac{m_0}{(1-p_T)}(T-p_T)s_T \right] \\
&=
\mathbb{E}\left[\frac{\one\{S=1\}}{\pi_S(\X,A=j)} \frac{p_A}{p(j)} (m_1-m_0-\tau_j) s_T\right] 
-\mathbb{E}\left[\frac{\one\{S=1\}}{\pi_S(\X,A=j)} \frac{p_A}{p(j)} \left\{\frac{T}{p_T}m_1-\frac{1-T}{1-p_T}m_0\right\} s_T\right],
\end{align*}
where the last holds from $\mathbb{E}[s_T|\X,\Z,S]=0$.
Moreover, we have
\begin{align*}
&\mathbb{E}[GH s_T]
=
\mathbb{E}[GE[H|\X,\Z,T,S] s_T]
=
\mathbb{E}\left[G
\left\{\frac{T}{p_T}m_1-\frac{1-T}{1-p_T}m_0\right\} s_T\right].
\end{align*}
Define
$$
\varphi_T = \frac{\one\{S=1\}}{\pi_S(\X,A=j)} \frac{p_A}{p(j)} (m_1-m_0-\tau_j).
$$
Combining above results, $\varphi_T$ clearly satisfies
$$
\mathbb{E}[\varphi_T s_Y] = \mathbb{E}[\varphi_T s_A] = \mathbb{E}[\varphi_T s_q] = 0.
$$

Next, we consider the contribution of $s_A$, which appears in $\mathbb{E}[{\dot G}H]$ and $\mathbb{E}[GHs_A]$.
Since it holds that
$$
\dot{p}_A = \mathbb{E}[(\one\{A=j\}-p_A)s_A|\X,\Z,S=1],
$$
we have
$$
\dot{G} = \frac{\partial G}{\partial p_A} \dot{p}_A
= \frac{\one\{S=1\}}{p_S} \frac{1}{p(j)} \mathbb{E}[(\one\{A=j\}-p_A)s_A|\X,\Z,S=1].
$$
Then, we have
\begin{align*}
\mathbb{E}[\dot{G}H] =& \mathbb{E}\left[\frac{\one\{S=1\}}{\pi_S} \frac{H}{p(j)} \mathbb{E}[(\one\{A=j\}-p_A)s_A|\X,\Z,S=1] \right] \\
=&
\mathbb{E}\left[\mathbb{E}\left[\frac{\one\{S=1\}}{\pi_S} \frac{H}{p(j)} \bigg|\X,\Z\right] (\one\{A=j\}-p_A)\frac{\one\{S=1\}}{\mathbb{P}(S=1|\X,\Z)}s_A \right]  \\
=&
\mathbb{E}\left[\frac{\one\{S=1\}}{\pi_S}\frac{\one\{A=j\}-p_A}{p(j)} (m_1-m_0-\tau_j)  s_A \right].
\end{align*}
From Assumption 4 (Area ignorability) and property of score
$$
\mathbb{E}[GHs_A] = \mathbb{E}[GHE[s_A|\X,\Z,S]]=0.
$$
Then, we define
$$
\varphi_A = \frac{\one\{S=1\}}{\pi_S}\frac{\one\{A=j\}-p_A}{p(j)} (m_1-m_0-\tau_j),
$$
which clearly satisfies 
$$
\mathbb{E}[\varphi_As_Y] = \mathbb{E}[\varphi_As_T] = \mathbb{E}[\varphi_As_q] = 0.
$$

Lastly,  we consider the contribution of $s_q$
Since $s_q$ is the score for $p(\x,\z,s)$ and satisfies $\mathbb{E}[f(\x,\z,s)s_q]=0$ for any square-integrable function $f$, we have 
$$
\mathbb{E}[GHs_q]=\mathbb{E}[\mathbb{E}[GH|\X.\Z,S]s_q]=0.
$$

Therefore, combining above results,
\begin{align*}
\varphi(\W)
=
&
\frac{\one\{S=1\}}{\pi_S} \frac{p_A}{p(j)}
\left\{
\frac{T(Y-m_1)}{p_T}
-
\frac{(1-T)(Y-m_0)}{1-p_T}
+m_1-m_0-\tau_j
\right\} \\
&+
\frac{\one\{S=1\}}{\pi_S}\frac{\one\{A=j\}-p_A}{p(j)} (m_1-m_0-\tau_j)
\end{align*}
is an efficient influence function.

\section{Proof of Theorem \ref{thm:an}}

At first, we give the additional assumptions for the proof.

\paragraph{Assumptions}

\begin{itemize}

\item[(SA1)] (Moments)
$\mathbb{E}[Y^2]<\infty$.

\item[(SA2)] (Rates)
$$
\|\widehat{m}_t - m_t\| = o_p(n^{-1/4}), \quad
\|\widehat{e} - e\| = o_p(n^{-1/4}), \quad
\|\widehat{\pi}_A - \pi_A\| = o_p(n^{-1/4}).
$$

\item[(SA3)] (Cross-fitting independence)
The nuisance estimator trained without observation $i$, $\widehat{\bta}^{(-k(i))}$ is independent of $\W_i$.
\end{itemize}

\paragraph{Proof}

For a measurable function $f$ and a random variable $X$, define $\mathbb{\mathbb{P}}_n f(X) = n^{-1} \sum_{i=1}^n f(X_i)$ and $\mathbb{\mathbb{P}} f(X) = \mathbb{E}[f(X_i)]$.
Horvitz--Thompson type estimator can be written as $\widehat{\tau}_j^{\rm HT} = \mathbb{\mathbb{P}}_n \phi^{(j)}(\W;\widehat{\bta})$, and can be decomposed as
$$
\widehat{\tau}_j^{\rm HT} - \tau_j
=
(\mathbb{\mathbb{P}}_n - \mathbb{\mathbb{P}})\phi^{(j)}(\W;\bta_0)
+
\mathbb{\mathbb{P}}_n[\phi^{(j)}(\W;\widehat{\bta})-\phi^{(j)}(\W;\bta_0)]
\equiv
(\mathbb{\mathbb{P}}_n - \mathbb{\mathbb{P}})\phi^{(j)}(\W;\bta_0)
+ 
R_n.
$$

Furthermore, $R_n$ can be decomposed as
$$
R_n
=
(\mathbb{\mathbb{P}}_n - \mathbb{\mathbb{P}})[\phi^{(j)}(\W;\widehat{\bta})-\phi^{(j)}(\W;\bta_0)]
+
\mathbb{\mathbb{P}}[\phi^{(j)}(\W;\widehat{\bta})-\phi^{(j)}(\W;\bta_0)]
\equiv
R_{n1}+R_{n2}.
$$

At first, we control the empirical process term $R_{n1}$.
From Assumption SA3 (Cross-fitting independence), conditional on the sigma-field generated by all cross-fitted
nuisance estimators, the summands are i.i.d. with mean zero.
Then, we have
$$
\mathbb{E}
\left[
\left\{
(\mathbb{\mathbb{P}}_n - \mathbb{\mathbb{P}})
[\phi^{(j)}(\W;\widehat{\bta}) - \phi^{(j)}(\W;\bta_0)]
\right\}^2
\right]
\le
\frac{1}{n}
\mathbb{E}
\left[
\left\{\phi^{(j)}(\W;\widehat{\bta}) - \phi^{(j)}(\W;\bta_0)\right\}^2
\right].
$$
By smoothness of $\phi^{(j)}$ and Assumption 3 (Overlap),
$$
\{\phi^{(j)}(\W;\widehat{\bta}) - \phi^{(j)}(\W;\bta_0)\}^2
\lesssim
\|\widehat{m}_0-m_0\|^2 + \|\widehat{m}_1-m_1\|^2 + \|\widehat{e}-e\|^2 + \|\widehat{\pi}_A-\pi_A\|^2.
$$
Using Assumption SA2 (Rates), we obtain
$$
\mathbb{E}
\left[
\{\phi^{(j)}(\W;\widehat{\bta}) - \phi^{(j)}(\W;\bta_0)\}^2
\right]
=
o_p(1),
$$
which implies $R_{n1}=o_p(n^{-1/2})$.

Next, we consider $R_{n2}$.
By a first-order Taylor expansion,
$$
R_{n2}
=
\partial_\bta \mathbb{\mathbb{P}}[\phi^{(j)}(\W;\bta_0)]
\cdot (\widehat{\bta}-\bta_0)
+
R_{n2}^{(2)},
$$
where $R_{n2}^{(2)}$ is a second-order remainder.
By the orthogonality property of $\phi^{(j)}$, namely $\partial_\bta \mathbb{\mathbb{P}}[\phi^{(j)}(\W;\bta_0)] = 0$, the first-order term vanishes and we obtain $R_{n2}=R_{n2}^{(2)}$.
Since the second-order remainder satisfies $R_{n2}^{(2)}= o_p(n^{-1/2})$ from Assumption SA2, we have $R_{n2}= o_p(n^{-1/2})$.

By combining above results, we have $R_n = o_p(n^{-1/2})$, and therefore
$$
\sqrt{n}(\widehat{\tau}_j^{\rm HT}-\tau_j)
=
\frac{1}{\sqrt{n}}\sum_{i=1}^n
(\phi^{(j)}(\W_i;\bta_0)-\tau_j)
+
o_p(1).
$$
Variance of $\widehat{\tau}_j^{\rm HT}$ can be written as
$$
\Var(\widehat{\tau}_j^{\rm HT})
=
\frac{1}{n}
\mathbb{E}[(\phi^{(j)}(\W_i;\bta_0)-\tau_j)^2]
+
o(n^{-1}).
$$

Next, we consider H\'ajek type estimator
$\widehat{\tau}_j^{\rm Hajek}
=
\mathbb{\mathbb{P}}_n \phi^{(j)}(\W;\widehat{\bta}) /
\mathbb{\mathbb{P}}_n w_{j1}(\W;\widehat{\bta})$.
From the similar arguments for HT type estimator, the numerator and denominator of $\widehat{\tau}_j^{\rm Hajek}$ can respectively be expanded as
$$
\mathbb{\mathbb{P}}_n \phi^{(j)}(\W;\widehat{\bta})
=
\mathbb{\mathbb{P}}\phi^{(j)}(\W;\bta_0)
+
(\mathbb{\mathbb{P}}_n-\mathbb{\mathbb{P}})\phi^{(j)}(\W;\bta_0)
+
o_p(n^{-1/2}), 
$$
and
$$
\mathbb{\mathbb{P}}_n w_{j1}(\W;\widehat{\bta})
=
\mathbb{\mathbb{P}}w_{j1}(\W;\bta_0)
+
(\mathbb{\mathbb{P}}_n-\mathbb{\mathbb{P}})w_{j1}(\W;\bta_0)
+
o_p(n^{-1/2}).
$$
Then, by using Taylor expansion,
\begin{align*}
&\widehat{\tau}_j^{\rm Hajek}
=
\frac{\mathbb{\mathbb{P}}_n \phi^{(j)}(\W;\widehat{\bta})}{\mathbb{\mathbb{P}}_n w_{j1}(\W;\widehat{\bta})}
=
\frac{\mathbb{\mathbb{P}}\phi^{(j)}(\W;\bta_0)}{\mathbb{\mathbb{P}}w_{j1}(\W;\bta_0)} 
+
\frac{1}{\mathbb{\mathbb{P}}w_{j1}(\W;\bta_0)}\{\mathbb{\mathbb{P}}_n \phi^{(j)}(\W;\widehat{\bta})-\mathbb{\mathbb{P}}\phi^{(j)}(\W;\bta_0)\} \\
&\qquad
-
\frac{\mathbb{\mathbb{P}}\phi^{(j)}(\W;\bta_0)}{\{\mathbb{\mathbb{P}}w_{j1}(\W;\bta_0)\}^2}\{\mathbb{\mathbb{P}}_n w_{j1}(\W;\widehat{\bta})-\mathbb{\mathbb{P}}w_{j1}(\W;\bta_0)\}
+
o_p(n^{-1/2}) \\
&=
\tau_j 
+
(\mathbb{\mathbb{P}}_n - \mathbb{\mathbb{P}}) \left[
\frac{\phi^{(j)}(\W;\bta_0)}{\mathbb{\mathbb{P}}w_{j1}(\W;\bta_0)}-\frac{\mathbb{\mathbb{P}}\phi^{(j)}(\W;\bta_0)}{\{\mathbb{\mathbb{P}}w_{j1}(\W;\bta_0)\}^2}w_{j1}(\W;\bta_0)
\right]
+o_p(n^{-1/2}).
\end{align*}
Therefore, we have
$$
\sqrt{n}(\widehat{\tau}_j^{\rm Hajek}-\tau_j)
=
\frac{1}{\sqrt{n}}
\sum_{i=1}^n
\frac{
\phi^{(j)}(\W;\bta_0)
-
\tau_j w_{j1}(\W_i;\bta_0)
}{\mathbb{E}[w_{j1}(\W_i;\bta_0)]}
+
o_p(1),
$$
and variance of $\widehat{\tau}_j^{\rm Hajek}$ can be written as
$$
\Var(\widehat{\tau}_j^{\rm Hajek})
=
\frac{1}{n}
\mathbb{E}
\left[
\left(
\frac{
\phi^{(j)}(\W;\bta_0)
-
\tau_j w_{j1}(\W_i;\bta_0)
}{\mathbb{E}[w_{j1}(\W_i;\bta_0)]}
\right)^2
\right]
+
o(n^{-1}).
$$

\section{Auxiliary results of Section \ref{sec:sim} and \ref{sec:rda}}

\subsection{Additional simulation results of Section \ref{sec:sim}}

\begin{figure}[h]
  \centering
  \begin{minipage}{0.4\columnwidth}
     \centering
     \includegraphics[width=\columnwidth]{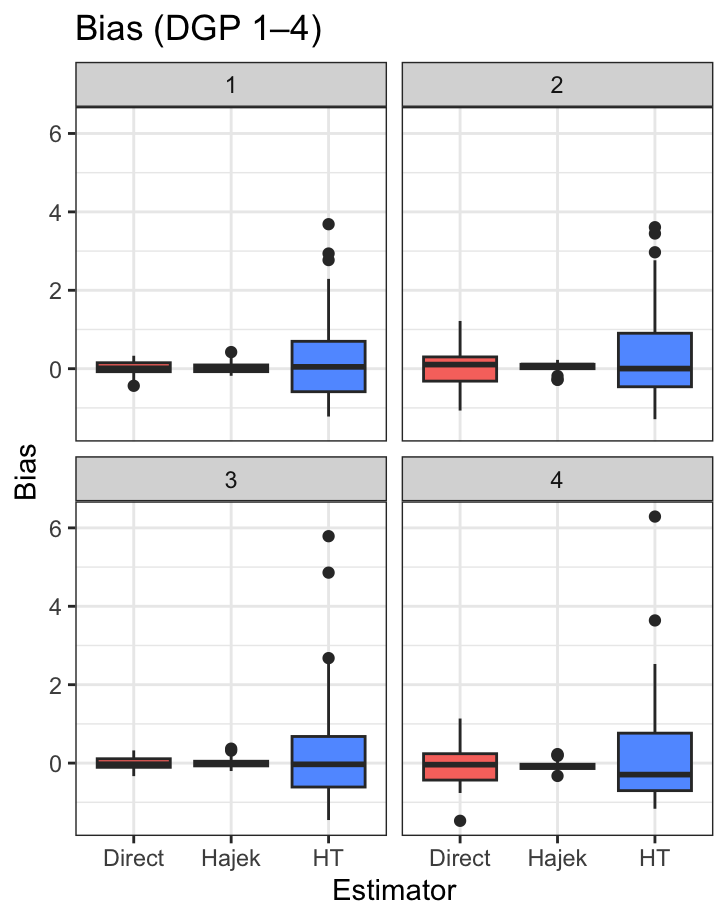}
  \end{minipage}
  \begin{minipage}{0.4\columnwidth}
     \centering
     \includegraphics[width=\columnwidth]{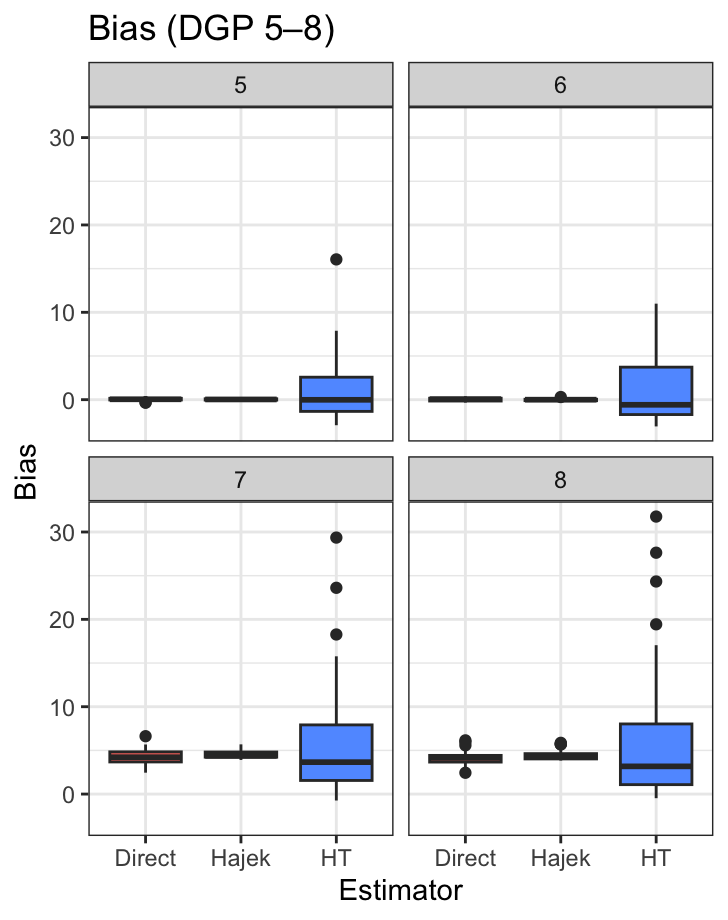}
  \end{minipage}
     \caption{Boxplots of the area-specific estimation bias across the 50 areas for each data generating process (DGP) with total population size $N=800{,}000$.  
For each estimator, the distribution of $\mathrm{Bias}(\hat{\tau}_j)$ across areas is summarized over the Monte Carlo
replications. }
     \label{fig:s4}
\end{figure}

\begin{figure}[h]
  \centering
  \begin{minipage}{0.4\columnwidth}
     \centering
     \includegraphics[width=\columnwidth]{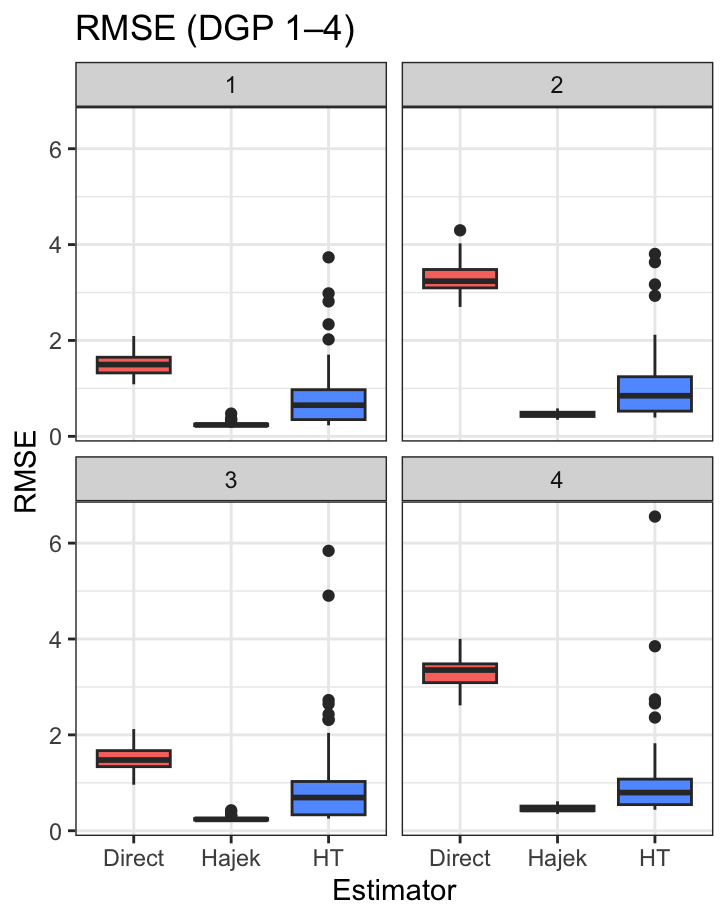}
  \end{minipage}
  \begin{minipage}{0.4\columnwidth}
     \centering
     \includegraphics[width=\columnwidth]{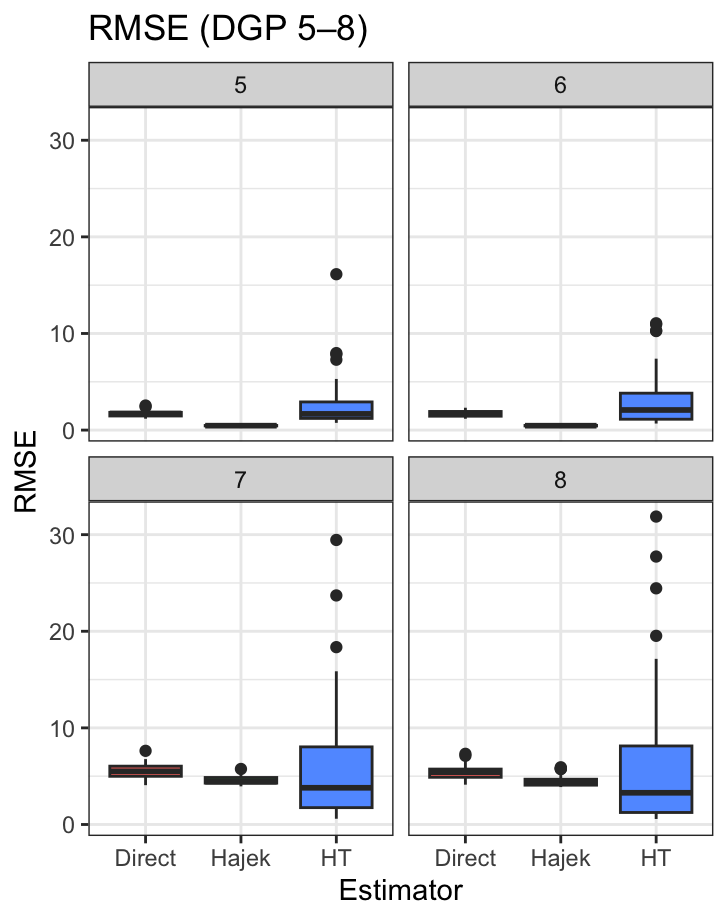}
  \end{minipage}
     \caption{Boxplots of the root mean squared error (RMSE) of the area-specific estimation across the 50 areas for each data generating process (DGP) with total population size $N=800{,}000$. 
For each estimator, the distribution of $\mathrm{RMSE}(\hat{\tau}_j)$ across areas is summarized over the Monte Carlo
replications. }
     \label{fig:s5}
\end{figure}

\begin{figure}[h]
  \centering
  \begin{minipage}{0.4\columnwidth}
     \centering
     \includegraphics[width=\columnwidth]{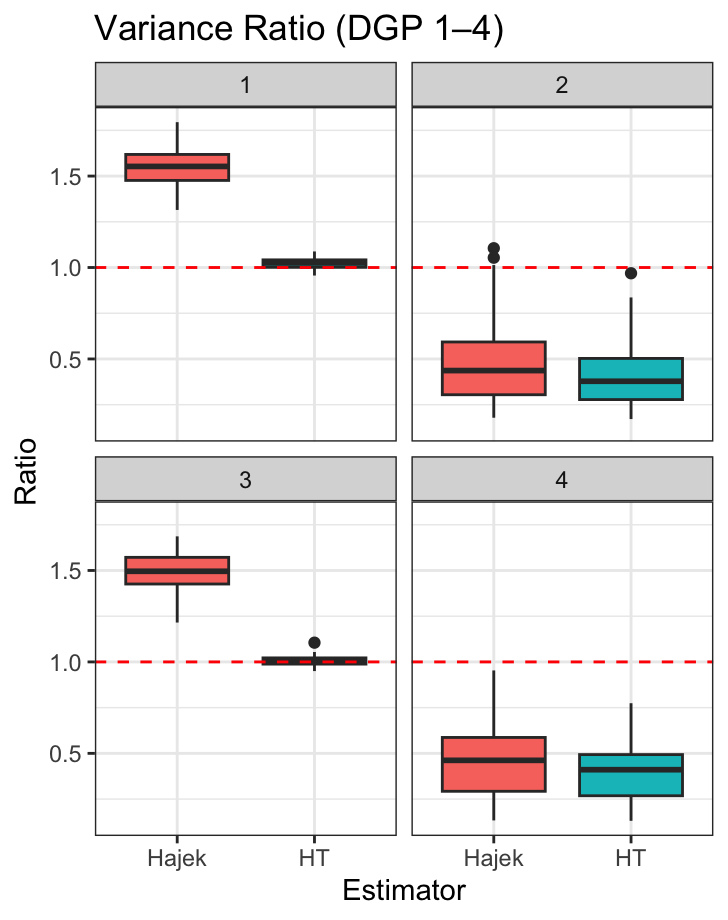}
  \end{minipage}
  \begin{minipage}{0.4\columnwidth}
     \centering
     \includegraphics[width=\columnwidth]{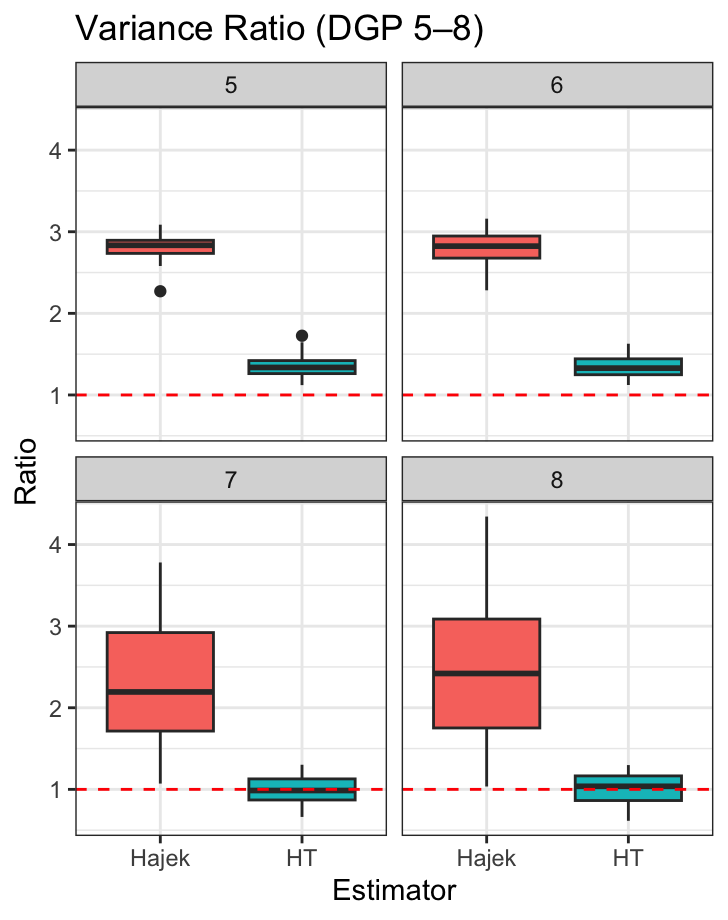}
  \end{minipage}
     \caption{Boxplots of the ratio of variance estimation of the area-specific estimation to the true variance across the 50 areas for each data generating process (DGP) with total population size $N=800{,}000$. 
For each estimator, the distribution of $\mathrm{Ratio}(\hat{V}_j)$ across areas is summarized over the Monte Carlo
replications. }
     \label{fig:s6}
\end{figure}

In this subsection, we give the performances of the proposed methods for the total population size $N=800{,}000$.  
Figures \ref{fig:s4}--\ref{fig:s6} respectively show boxplots of the bias, RMSE, and variance estimation accuracy across the 50 areas for each data generating process, and Table \ref{tab:80} summarizes the corresponding averages and standard deviations.

\begin{table}[h]
\caption{Average value of Bias, RMSE and Variance Ratio for each data generating process (DGP) with total population size $N=800{,}000$. Values in parentheses are standard deviations.}
\centering

\medskip
\begin{tabular}{
c
S[table-format=1.3] S[table-format=1.3] S[table-format=1.3]
S[table-format=1.3] S[table-format=1.3] S[table-format=1.3]
S[table-format=1.3] S[table-format=1.3]
}
\toprule
\multirow{2}{*}{DGP} & \multicolumn{3}{c}{Bias} & \multicolumn{3}{c}{RMSE} & \multicolumn{2}{c}{Variance Ratio} \\
\cmidrule(lr){2-4} \cmidrule(lr){5-7} \cmidrule(lr){8-9}
& {Direct} & {HT} & {H\'ajek} & {Direct} & {HT} & {H\'ajek} & {HT} & {H\'ajek} \\
\midrule

\multirow{2}{*}{DGP1}
& 0.019 & 0.275 & 0.020 & 1.485 & 0.864 & {\bfseries 0.250} & 1.023 & 1.543 \\
& \multicolumn{1}{c}{\footnotesize (0.174)} & \multicolumn{1}{c}{\footnotesize (1.085)} & \multicolumn{1}{c}{\footnotesize (0.125)}
& \multicolumn{1}{c}{\footnotesize (0.229)} & \multicolumn{1}{c}{\footnotesize (0.765)} & \multicolumn{1}{c}{\footnotesize (0.046)}
& \multicolumn{1}{c}{\footnotesize (0.032)} & \multicolumn{1}{c}{\footnotesize (0.102)} \\
\addlinespace[3pt]

\multirow{2}{*}{DGP2}
& 0.052 & 0.348 & 0.052 & 3.284 & 1.107 & {\bfseries 0.462} & 0.419 & 0.487 \\
& \multicolumn{1}{c}{\footnotesize (0.493)} & \multicolumn{1}{c}{\footnotesize (1.208)} & \multicolumn{1}{c}{\footnotesize (0.116)}
& \multicolumn{1}{c}{\footnotesize (0.333)} & \multicolumn{1}{c}{\footnotesize (0.822)} & \multicolumn{1}{c}{\footnotesize (0.060)}
& \multicolumn{1}{c}{\footnotesize (0.188)} & \multicolumn{1}{c}{\footnotesize (0.239)} \\
\addlinespace[3pt]

\multirow{2}{*}{DGP3}
& 0.005 & 0.325 & 0.005 & 1.498 & 1.031 & {\bfseries 0.249} & 1.006 & 1.486 \\
& \multicolumn{1}{c}{\footnotesize (0.163)} & \multicolumn{1}{c}{\footnotesize (1.453)} & \multicolumn{1}{c}{\footnotesize (0.115)}
& \multicolumn{1}{c}{\footnotesize (0.258)} & \multicolumn{1}{c}{\footnotesize (1.112)} & \multicolumn{1}{c}{\footnotesize (0.043)}
& \multicolumn{1}{c}{\footnotesize (0.028)} & \multicolumn{1}{c}{\footnotesize (0.115)} \\
\addlinespace[3pt]

\multirow{2}{*}{DGP4}
& -0.078 & 0.196 & -0.078 & 3.305 & 1.085 & {\bfseries 0.462} & 0.395 & 0.453 \\
& \multicolumn{1}{c}{\footnotesize (0.472)} & \multicolumn{1}{c}{\footnotesize (1.368)} & \multicolumn{1}{c}{\footnotesize (0.114)}
& \multicolumn{1}{c}{\footnotesize (0.315)} & \multicolumn{1}{c}{\footnotesize (1.039)} & \multicolumn{1}{c}{\footnotesize (0.061)}
& \multicolumn{1}{c}{\footnotesize (0.139)} & \multicolumn{1}{c}{\footnotesize (0.177)} \\
\addlinespace[3pt]

\multirow{2}{*}{DGP5}
& 0.053 & 0.025 & 0.024 & 1.698 & 2.726 & {\bfseries 0.456} & 1.348 & 2.814 \\
& \multicolumn{1}{c}{\footnotesize (0.170)} & \multicolumn{1}{c}{\footnotesize (3.606)} & \multicolumn{1}{c}{\footnotesize (0.108)}
& \multicolumn{1}{c}{\footnotesize (0.305)} & \multicolumn{1}{c}{\footnotesize (2.672)} & \multicolumn{1}{c}{\footnotesize (0.036)}
& \multicolumn{1}{c}{\footnotesize (0.124)} & \multicolumn{1}{c}{\footnotesize (0.139)} \\
\addlinespace[3pt]

\multirow{2}{*}{DGP6}
& 0.026 & -0.015 & -0.044 & 1.711 & 3.012 & {\bfseries 0.462} & 1.344 & 2.808 \\
& \multicolumn{1}{c}{\footnotesize (0.191)} & \multicolumn{1}{c}{\footnotesize (3.853)} & \multicolumn{1}{c}{\footnotesize (0.137)}
& \multicolumn{1}{c}{\footnotesize (0.299)} & \multicolumn{1}{c}{\footnotesize (2.705)} & \multicolumn{1}{c}{\footnotesize (0.046)}
& \multicolumn{1}{c}{\footnotesize (0.124)} & \multicolumn{1}{c}{\footnotesize (0.181)} \\
\addlinespace[3pt]

\multirow{2}{*}{DGP7}
& 4.252 & 4.526 & 4.473 & 5.523 & 6.166 & {\bfseries 4.576} & 0.990 & 2.354 \\
& \multicolumn{1}{c}{\footnotesize (0.820)} & \multicolumn{1}{c}{\footnotesize (6.351)} & \multicolumn{1}{c}{\footnotesize (0.421)}
& \multicolumn{1}{c}{\footnotesize (0.774)} & \multicolumn{1}{c}{\footnotesize (6.257)} & \multicolumn{1}{c}{\footnotesize (0.424)}
& \multicolumn{1}{c}{\footnotesize (0.170)} & \multicolumn{1}{c}{\footnotesize (0.747)} \\
\addlinespace[3pt]

\multirow{2}{*}{DGP8}
& 4.151 & 4.427 & 4.657 & 5.407 & 6.287 & {\bfseries 4.475} & 1.010 & 2.471 \\
& \multicolumn{1}{c}{\footnotesize (0.784)} & \multicolumn{1}{c}{\footnotesize (7.449)} & \multicolumn{1}{c}{\footnotesize (0.518)}
& \multicolumn{1}{c}{\footnotesize (0.725)} & \multicolumn{1}{c}{\footnotesize (7.365)} & \multicolumn{1}{c}{\footnotesize (0.521)}
& \multicolumn{1}{c}{\footnotesize (0.193)} & \multicolumn{1}{c}{\footnotesize (0.895)} \\

\bottomrule
\end{tabular}
\label{tab:80}
\end{table}

\subsection{Details of variables in Section \ref{sec:rda}}

The outcome variable is based on respondents' feeling thermometer ratings toward the two major presidential candidates in the 2024 election.
The respondents evaluated each candidate on a 0--100 scale.
We construct the outcome as the difference between the assessment of Donald Trump and that of Kamala Harris, $Y_i = FT_i^{Trump} - FT_i^{Harris}$, which takes values between $-100$ and $100$.
Positive values indicate warmer feelings toward Trump relative to Harris, while negative values indicate warmer feelings toward Harris.

The treatment variable $T_i$ indicates whether the respondent was contacted by a political campaign during the election cycle.
Specifically, we define treatment as having received campaign contact from both the Democratic Party and the Republican Party.
Our goal is to estimate the average treatment effect of campaign contact on the outcome variable separately for each of the 51 U.S. states including the District of Columbia.

We consider two sets of covariates.
The first set, denoted by $\X_i$, consists of demographic variables that are available both in the survey and from external population data sources such as the U.S. Census. 
These variables include age, sex (1. Male, 0. Female), race and ethnicity (1. White, non-Hispanic, 2. Black, non-Hispanic, 3. Hispanic, 4. Asian or Native Hawaiian/Other Pacific Islander, non-Hispanic, 5. Native American/Alaska Native or other race, non-Hispanic, 6. Multiple races, non-Hispanic), educational attainment (1. High school or less, 2. Some college or associate degree, 3. Bachelor's degree, 4. Graduate degree) and total household income (1. Under \$9{,}999, 2. \$10{,}000--\$29{,}999, 3. \$30{,}000--\$59{,}999, 4. \$60{,}000--\$99{,}999, 5. \$100{,}000--\$249{,}999, 6. \$250{,}000 or more).

The second set of covariates, $\Z_i$, includes seven survey-based measures of political predispositions and behavior not available in census data: party identification, ideological self-placement, attention to politics, interest in following political campaigns, political knowledge, past political participation and trust in news media.
The identification of the party and the ideological self-placement are measured on seven-point scales (“1. Strong Democrat”–“7. Strong Republican”) and (“1. Extremely liberal”–“7. Extremely conservative”), respectively.
Attention to politics and interest in the following campaigns are measured on five- and three-point scales, respectively, with higher values indicating greater political interest. 
Political knowledge is defined as the standardized sum of four binary indicators for correct responses to factual questions about U.S. politics.
Past political participation is captured by an indicator for voting in the 2020 presidential election (1 if voted, 0 otherwise). 
Trust in the news media is measured on a five-point scale, with higher values indicating greater trust.
Summary statistics are reported in Table \ref{tab:sum}.

\begin{table}[htbp]
\centering
\caption{Summary Statistics by Treatment Status}
\begin{tabular}{lcccccc}
\toprule
& \multicolumn{3}{c}{$T_i=0$} & \multicolumn{3}{c}{$T_i=1$} \\
\multirow{2}{*}{Variable} & \multirow{2}{*}{Mean} & \multicolumn{2}{c}{SD} & \multirow{2}{*}{Mean} & \multicolumn{2}{c}{SD}  \\
& & (overall) & (between) & & (overall) & (between) \\
\midrule
Feeling thermometer & -2.84 & 75.0 & 27.0 & -8.85 & 70.6 & 35.9 \\
Age & 54.3 & 16.8 & 4.6 & 50.9 &16.8 & 8.1 \\
Sex & 0.53 & 0.50 & 0.13 & 1.49 & 0.50 & 0.31 \\
Race & 1.33 & 1.01 & 0.24 & 1.34 & 1.07 & 0.68 \\
Education & 2.58 & 1.02 & 0.37 & 2.79 & 0.99 & 0.57 \\
Income & 4.01 & 1.32 & 0.39 & 4.29 & 1.27 & 0.58 \\
Party identification & 4.01 & 2.33 & 0.92 & 3.84 & 2.29 & 1.20 \\
Ideological self-placement & 4.16 & 1.71 & 0.74 & 4.08 & 1.61 & 0.86 \\
Attention to politics & 3.67 & 1.02 & 0.37 & 3.68 & 0.98 & 0.65 \\
Interest in campaigns & 2.44 & 0.68 & 0.27 & 2.43 & 0.69 & 0.34 \\
Political knowledge & 0.18 & 0.93 & 0.24 & 0.20 & 0.93 & 0.51 \\
Past political participation & 0.84 & 0.36 & 0.09 & 0.89 & 0.31 & 0.23 \\
Trust in news media & 3.59 & 1.08 & 0.27 & 3.55 & 1.07 & 0.61 \\
\bottomrule
\end{tabular}
\label{tab:sum}
\end{table}

\end{document}